\documentclass[12pt]{amsproc}

\usepackage{amsmath}
\usepackage{amssymb}
\usepackage{amsthm}

\newtheorem{theorem}{Theorem}[section]
\newtheorem{lemma}[theorem]{Lemma}

\newtheorem{proposition}[theorem]{Proposition}
\newtheorem{corollary}[theorem]{Corollary}
\theoremstyle{definition}
\newtheorem{definition}[theorem]{Definition}

\theoremstyle{remark}
\newtheorem{remark}[theorem]{Remark}

\numberwithin{equation}{section}

\newcommand{\ord}{{\mathrm{ord}}}
\newcommand{\Sol}{{\mathrm{Sol}}}
\newcommand{\Res}{\operatornamewithlimits{Res}}
\newcommand{\C}{{\mathbb C}}
\newcommand{\Z}{{\mathbb Z}}
\newcommand{\R}{{\mathbb R}}
\newcommand{\J}{{\mathcal J}}
\newcommand{\N}{{\mathbb N}}
\newcommand{\selh}{\hat{\mathfrak{sl}}(2)}
\newcommand{\sel}{\mathfrak{sl}(2)}

\begin{document}

\title[Bethe Equation at $q=0$]{Bethe Equation at $q=0$, M\"obius Inversion
Formula, and Weight Multiplicities:\\
 I.\ $\sel$ case
\footnote{
For {\em Proceedings of Workshop on
Physical Combinatorics}, IIAS, Kyoto, 1999.}
}\
\author{Atsuo Kuniba}
\address{Institute of Physics, University of Tokyo, Tokyo 153-8902, Japan}
\email{atsuo@hep1.c.u-tokyo.ac.jp}
%
\author{Tomoki Nakanishi}
\address{Graduate School of Mathematics,
    Nagoya University, Nagoya 464-8602, Japan}
\email{nakanisi@math.nagoya-u.ac.jp}
%
%

\begin{abstract}
The $U_q(\selh)$ Bethe equation is studied at $q=0$.
A linear congruence equation is proposed related to the 
string solutions.
The number of its off-diagonal solutions is 
expressed in terms of an explicit combinatorial formula and 
coincides with the weight multiplicities of the quantum space.
\end{abstract}

\maketitle

\section{Introduction}\label{sec:intro}

\subsection{Background}\label{subsec:back}

Consider the periodic spin $\frac{1}{2}$ XXX Heisenberg hamiltonian
\begin{equation*}
H_{XXX} = J \sum_{n=1}^L(\sigma^x_n\sigma^x_{n+1} + 
\sigma^y_n\sigma^y_{n+1} + \sigma^z_n\sigma^z_{n+1})
\end{equation*}
acting on the tensor product of the $L$-copies of the 
vector representations of ${\mathfrak{sl}}(2)$:
\begin{equation*}
W = \C^2 \otimes \cdots \otimes \C^2.
\end{equation*}
Since $H_{XXX}$ is ${\mathfrak{sl}}(2)$-linear, its spectrum is 
degenerated within  the irreducible components in the decomposition:
\begin{equation*}
W = \bigoplus_{\lambda \in (\Z_{\ge 0})\Lambda_1} [W: V_\lambda] \; V_\lambda,
\end{equation*}
where $\Lambda_1$ is the fundamental weight and 
$V_\lambda$ denotes the irreducible module with highest weight $\lambda$.
Diagonalization of $H_{XXX}$ was achieved by Bethe \cite{Be} in 1931.
Associated to each solution of the simultaneous equations 
(Bethe equation)
\begin{equation*}
\left(\frac{u_j + \sqrt{-1}}{u_j - \sqrt{-1}}\right)^L
= -\prod_{j=1}^M 
\frac{u_j-u_k+2\sqrt{-1}}{u_j-u_k-2\sqrt{-1}}
\qquad j = 1, \ldots, M,
\end{equation*}
he proposed (Bethe ansatz) a vector $\psi \in W$ (Bethe vector) such that
\begin{align*}
H_{XXX} \psi &= E \psi \qquad E \in \C,\\
\left(\sum_{n=1}^L\sigma^+_n\right) \psi &= 0 \qquad 
\sigma^+_n = \sigma^x_n + \sqrt{-1}\sigma^y_n,\\
\left(\sum_{n=1}^L\sigma^z_n\right) \psi &= (L-2M)\psi \qquad 0 \le M \le
[\frac{L}{2}].
\end{align*}
The second and the third properties  (cf.\ \cite{TF})
tell that the Bethe vector is ${\mathfrak{sl}}(2)$-highest of 
weight $(L-2M)\Lambda_1$.
Therefore in order to have the completeness of the Bethe ansatz,
there should exist as many solutions to the Bethe equation 
as the multiplicity of  $V_{(L-2M)\Lambda_1}$ in $W$,
$[W: V_{(L-2M)\Lambda_1}]  (=\binom{L}{M}-\binom{L}{M-1})$.
Actually $\psi$ can be vanishing depending on the solutions 
$\{u_1, \ldots, u_M\}$.
In particular, it is so if  $u_i = u_j$ for some $i \neq j$.

It was Bethe himself who studied the completeness 
with the introduction of strings. (He called it ``WellenKomplex".)
It is a solution in which $\{u_1, \ldots, u_M\}$ are arranged as 
\begin{equation}\label{eq:string}
\bigcup_{m\in \N}\bigcup_{1\le \alpha\le N_m}\bigcup_{u_{m\alpha} \in \R}
\{ u_{m \alpha} + \sqrt{-1}(m+1-2i) + \epsilon_{m\alpha i}
\mid 1 \le i \le m \}
\end{equation}
for each partition $M = \sum_{m \in \N}mN_m\; (N_m \in \Z_{\ge 0})$.
Here $\N = \Z_{\ge 1}$ denotes the set of 
positive integers and $\epsilon_{m\alpha i}$ stands for a small deviation.
The $m$-tuple configuration (with negligible $\epsilon_{m\alpha i}$) 
is called the $m$-string with string center $u_{m\alpha}$.
$N_m$ is the number of $m$-strings.
In general to expect such a behavior for the solutions 
is called the string hypothesis.
Actually in a strict sense, it is known invalid as exemplified 
already for $M = 2$ and $L > 21$ (cf.\ \cite{EKS,JD}).
Nevertheless Bethe's counting of the number of string solutions led to the 
discovery of the identity ($M \le [\frac{L}{2}]$):
\begin{equation}\label{eq:completeness0}
\sum_N \prod_{m \in \N}
\binom{L-2\sum_{k\ge 1}\min(m,k)N_k + N_m}{N_m} = 
[W: V_{(L-2M)\Lambda_1}],
\end{equation}
where $\sum_N$ runs over $N_1$, $N_2$, $\ldots \in \Z_{\ge 0}$
such that $M = \sum_{m \ge 1}mN_m$.
In his counting each summand in the left side 
represents the number of string solutions 
corresponding to the prescribed values of 
$N_1, N_2, \ldots$.
The binomial coefficients are originated in the fermionic restriction
on the solutions $u_i \neq u_j\, (i \neq j)$.
The expression like the left side is called the fermionic formula and 
the above identity is called the combinatorial completeness
of the string hypothesis.

Despite the gap with the completeness in the rigorous sense,
the above result opened a fruitful link between quantum integrable systems
and representation theory.
For a class of Bethe ansatz solvable models with Yangian symmetry 
$Y(X_n)$, one can set up fermionic formulae following Bethe's 
counting.
If the combinatorial completeness holds, they yield the 
multiplicities of irreducible $X_n$-modules in
tensor products of a variety of finite dimensional 
$Y(X_n)$-modules.
The XXX chain corresponds to the $Y({\mathfrak{sl}}(2))$ case.
The fermionic formula associated with $Y(X_n)$ in such a sense 
was firstly written down in \cite{KR} for general $X_n$, 
where the combinatorial completeness was also announced for 
the classical types 
$X_n = A_n$, $B_n$, $C_n$ and $D_n$.
The proof of the combinatorial completeness boils down to 
showing recursion relations ($Q$-system) among classical characters of 
certain $Y(X_n)$-modules (cf.\ \cite{HKOTY}).

\subsection{Present work}\label{subsec:present}

The XXX chain admits an integrable $q$-de\-for\-mation
 called the XXZ chain:
\begin{equation*}
H_{XXZ} = J \sum_{n=1}^L(\sigma^x_n\sigma^x_{n+1} + 
\sigma^y_n\sigma^y_{n+1} + \frac{q+q^{-1}}{2}\sigma^z_n\sigma^z_{n+1}).
\end{equation*}
In place of the Yangian $Y({\mathfrak{sl}}(2))$, 
the underlying symmetry of the model is the quantum affine 
algebra $U_q(\selh)$ as is well known.
Accordingly we regard the space $W$ (called the quantum space) as a 
$U_q(\selh)$-module.
Under the periodic boundary condition the spectrum is determined 
{}from the solutions of the Bethe equation
($i = 1$, \dots, $M$)
\begin{equation*}
\left(\frac{\sin{\pi\!\left(u_i + \sqrt{-1}\hbar\right)}}
{\sin{\pi\!\left(u_i - \sqrt{-1}\hbar\right)}}\right)^{L}
= - \prod_{j=1}^{M}
\frac{\sin\pi\!\left(u_i - u_j + 2\sqrt{-1}\hbar \right)}
{\sin\pi\!\left(u_i - u_j - 2\sqrt{-1}\hbar \right)},
\end{equation*}
where $\hbar$ is related to $q$ by $q = e^{-2\pi \hbar}$.
When the deformation parameter $q$ tends to $1$, 
the above equation reduces to the one in Section \ref{subsec:back}
by replacing $u_j$ by $\hbar u_j$ and setting $\hbar \rightarrow 0$.
A significant difference from the $q=1$ case is that 
the hamiltonian is no longer invariant under  ${\mathfrak{sl}}(2)$
nor $U_q({\mathfrak{sl}}(2))$
as far as a finite chain ($L<\infty$) is considered 
under the periodic boundary condition.
For the completeness, the number of solutions to the 
Bethe equation should therefore coincide with the weight 
multiplicity of $(L-2M)\Lambda_1$, 
$\dim W_{(L-2M)\Lambda_1} (=\binom{L}{M})$,
rather than  the multiplicity of $V_{(L-2M)\Lambda_1}$.
Similar facts are valid also for the generalized model 
in which  $W$ is replaced with 
\begin{equation*}
W(\nu) =\bigotimes_{s\geq 1}
(W_s)^{\otimes \nu_s},
\end{equation*}
where $\nu_s \in \Z_{\ge 0}$ and 
$W_s$ stands for the $(s+1)$-dimensional irreducible module.
See (\ref{eq:be}) for the corresponding Bethe equation.

The purpose of this paper is to study the 
Bethe equation and to formulate another version of the combinatorial 
completeness at $q=0$.
This is inspired by the crystal theory, where the 
simplification at $q=0$ is known to lead to many fascinating features.
In terms of the exponential variables $x_j = e^{2\pi\sqrt{-1}u_j}$,
we shall consider a class of meromorphic solutions $x_j = x_j(q)$ around 
$q=0$ which correspond to the strings.
In a sense we are approaching to the point $q=0$ 
within the off-critical regime $\hbar \in \R_{>0}$,
avoiding the parity and the arithmetic complexity of 
strings \cite{TS} in the critical regime $\hbar \in \sqrt{-1}\R$.

It is a routine calculation to reduce the Bethe equation to the one
for string centers for general $\nu = (\nu_s)$ and $N = (N_m)$.
At $q=0$ the resulting string center equation (SCE) is  
a linear congruence equation (\ref{eq:sce2}).
As a remnant of the fermionic restriction we seek the 
off-diagonal solutions ({\sc Definition} \ref{def:off2}) of them.
They are counted systematically by means of the M\"obius inversion formula.
When $P_m := \sum_{k \ge 1}\min(m,k)(\nu_k - 2N_k) \ge 0$ for any $m$ 
such that $N_m > 0$, we find that the result is expressed as 
(cf.\ {\sc Theorem} \ref{th:Rnds}):
\begin{align*}
R(\nu,N) &= \sum_{J \subset \N} D_J
\prod_{m \in \N \setminus J}\binom{P_m + N_m}{N_m}
\prod_{m \in J}\binom{P_m + N_m - 1}{N_m - 1},\\
D_J &= \begin{cases}
1 & \text{ if } J = \emptyset\\
\det_{m,k \in J}(2\min(m,k)-\delta_{m,k}) & \text{otherwise}.
\end{cases}
\end{align*}
In the XXZ case $\nu_s = L\delta_{s,1}$, the $J = \emptyset$ term here is 
equal to the summand in the 
left side of (\ref{eq:completeness0}).
With this $R(\nu,N)$ 
the combinatorial completeness at $q=0$ is stated  
as (cf.\ {\sc Theorem} \ref{th:completeness3}) 
\begin{equation*}
\sum_{N}{}^{(\lambda)} R(\nu,N) = \dim
W(\nu)_\lambda
\qquad \lambda \in
\Z\Lambda_1,
\end{equation*}
where the sum $\sum_{N}{}^{(\lambda)}$ runs over 
$N_1$, $N_2$, $\ldots \in \Z_{\ge 0}$ such that 
$\sum_{j \ge 1}j(\nu_j - 2N_j)\Lambda_1 = \lambda$.
This is an non-trivial identity even
when $\dim W(\nu)_\lambda = 0$ 
for $\lambda <0$.
Curiously, the left side in general involves the contributions from those
$N$ that 
break the $P_m \ge 0$ condition said above.
In the course of the proof we will clarify the relation 
between the known fermionic formula as in (\ref{eq:completeness0}) 
and our $R(\nu,N)$.
It is most transparently presented in terms of generating functions.
See (\ref{eq:RKKinf}) and (\ref{eq:RKK2}).

The layout of the paper is as follows.
In Section \ref{sec:beq0}, we study the Bethe equation at $q=0$.
We explain the relation between solutions of SCE and string solutions of the 
Bethe equation.
In Section \ref{sec:counting}
we derive the formula $R(\nu,N)$ by counting the off-diagonal solutions of SCE.
In Section \ref{sec:multiplicity} we prove the 
combinatorial completeness.
In Section \ref{sec:discussion} we give a summary and 
discussion.
Appendix \ref{app:mobius} is a summary of the M\"obius inversion trick 
used in Section \ref{sec:counting}.

A few remarks are in order.
In \cite{KL} combinatorial completeness was 
investigated when $q$ is a root of unity.
The fermionic formula there for weight multiplicities 
is different from ours.
We expect that their result describes the rich singular behavior of the 
meromorphic solutions (around $q=0$) of the Bethe equation on the 
convergence circle $\vert q \vert =1$.
In \cite{LS} SCE has been obtained for the XXZ case.
There is  a statement similar to the combinatorial completeness at $q=0$ 
without an explicit formula as $R(\nu,N)$.
In \cite{TV} they studied a deformation of the XXZ type Bethe equation
and showed the completeness of the Bethe vectors for the 
admissible off-diagonal solutions at a generic value of the 
deformation parameter.

Many results in this paper admit generalizations to 
$U_q(X^{(1)}_n)$ case for arbitrary $X_n$.
It will be a subject of our forthcoming paper.

\section{Bethe equation at $q=0$}\label{sec:beq0}

In this section  we start from the Bethe equation
and seek string solutions in
$q \rightarrow 0$ limit.
We introduce an equation for string centers  (SCE), which is
a linear congruence equation.
Later sections will be devoted to studies of the  SCE.
Our aim here is to explain the precise relation between
solutions of SCE and string solutions of the Bethe equation.
Our theorems mostly concern what we call generic string solutions.

We let $\alpha_1$ and $\Lambda_1$  denote the simple root and the
fundamental weight of $\sel$.
$\alpha_1 = 2\Lambda_1$.
In this paper $U_q(\selh)$ means the quantum affine algebra
without the derivation operator.
It is denoted by $U'_q(\selh)$ in some literature.

For a meromorphic function $f(q)$ around $q=0$,
we will use the notation
\begin{align*}
f(q) &= q^{\ord(f)}(f^0 + f^1 q + \cdots )
 \qquad \ord(f) \in \Z,\ f^0\neq 0,\\
%
%
\tilde{f}(q) &= q^{-\ord(f)}f(q) = f^0 + f^1q + \cdots.
\end{align*}
We call $\ord(f)$ the order, $f^0$ the leading coefficient of $f$.
$\tilde{f}$ is a holomorphic function around $q=0$
with nonzero constant term.
Note that $f^0 = \tilde{f}^0$.

\subsection{Bethe equation}

Consider a solvable vertex model associated with $U_q(\selh)$.
Let
\begin{equation}\label{eq:qspace}
\displaystyle W(\nu) =
\bigotimes_{s\geq 1}
(W_s)^{\otimes \nu_s}
\end{equation}
be the $U_q(\selh)$-module (the quantum space) on which the
commuting family of row-to-row
transfer matrices act.
We assume that only finitely many $\nu_s$'s are nonzero
throughout this paper.
Here $W_s$ stands for an $(s+1)$-dimensional irreducible module with 
highest weight $s\Lambda_1$.
Each $W_s$ depends on a spectral parameter, which may be interpreted as
inhomogeneity of the interaction.
The Bethe equation relevant to the spectrum of the
transfer matrices
also depends on those spectral parameters.
In this paper we concentrate on the regime and the situation in which
the Bethe equation takes the form ($i=1$, \dots, $M$):
\begin{equation}\label{eq:be}
\prod_{s \ge 1} \left(
\frac{\sin\pi\!\left(u_i + \sqrt{-1} s \hbar\right)}
{\sin\pi\!\left(u_i - \sqrt{-1} s\hbar\right)}\right)^{\nu_s}
= - \prod_{j=1}^{M}
\frac{\sin\pi\!\left(u_i - u_j + 2\sqrt{-1}\hbar \right)}
{\sin\pi\!\left(u_i - u_j - 2\sqrt{-1}\hbar \right)},
\end{equation}
Here $\hbar \in \R_{>0}$ and $M \in \Z_{\ge 0}$.
This is a regime in which the so-called parity \cite{TS}
is irrelevant.
Integer shifts of $u_j$ do not lead to
a new Bethe vector, hence one should consider $u_j \in \R/\Z$.
Setting
\begin{equation*}
q = e^{-2\pi\hbar},\quad  x_j=e^{2\pi\sqrt{-1}u_j},
\end{equation*}
(\ref{eq:be}) can be written as polynomial equations on $x_j$'s:
\begin{equation}\label{eq:ba0}
F_{i+}G_{i-}=-
F_{i-}G_{i+}
\qquad i=1, \dots, M,
\end{equation}
where
\begin{alignat*}{2}
F_{i+}&=
\prod_{s\geq 1}(x_iq^{{s}}
-1)^{\nu_s},
&\quad
G_{i+}&=
\prod_{j=1}^{M}
(x_iq^{2}
 -x_j),
\\
F_{i-}&=
\prod_{s\geq 1}(x_i
-q^{s})^{\nu_s},
&\quad
G_{i-}&=
\prod_{j=1}^{M}
(x_i
 -x_jq^{2}).
\end{alignat*}
The equation is invariant under the permutation
of the variables $x_i\leftrightarrow x_j$.

We are interested in meromorphic solutions
$(x_i)$, $x_i=x_i(q)$ of (\ref{eq:ba0}) around
$q=0$.
We set $x_i(q)=q^{d_i}z_i(q)$, where $d_i = \ord(x_i)$,
and $z_i(q) = \tilde{x}_i(q)$.
Then the Bethe equation for $z_i(q)$ is given by (\ref{eq:ba0})
with $F_{i \pm},\, G_{i\pm}$ now specified as
\begin{alignat}{2}
F_{i+}&=
\prod_{s\geq 1}(z_iq^{d_i+{s}}
-1)^{\nu_s},&\quad
G_{i+}&=
\prod_{j=1}^{M}
(z_iq^{d_i+2}
 -z_jq^{d_j}),
\\
F_{i-}&=
\prod_{s\geq 1}(z_iq^{d_i}
-q^{s})^{\nu_s},
&\quad
G_{i-}&=
\prod_{j=1}^{M}
(z_iq^{d_i}
 -z_jq^{d_j+2}).
\end{alignat}
This equation is invariant under the permutation
of the variables $z_i\leftrightarrow z_j$,
only when $d_i=d_j$.

\subsection{String Solution}

\begin{definition}\label{def:admissible}\upshape
A meromorphic solution $(x_i)$
of (\ref{eq:ba0}) is called
 {\em inadmissible (admissible)\/} if
$F_{i+}G_{i-} = F_{i-}G_{i+}=0$ for some $i$
as a function of $q$  around $q=0$ (otherwise).
\end{definition}
Let $N = (N_m)$ be an infinite sequence of non-negative integers such that
\begin{equation}\label{eq:mpartition}
M = \sum_{m \ge 1} mN_m.
\end{equation}
\begin{definition}\label{def:string}\upshape
A meromorphic solution $(x_i)$ of (\ref{eq:ba0})
is called a
{\em string solution of pattern
$N=(N_m)$} if \hfill\break
(i) $(x_i)$ is admissible.\hfill\break
(ii) $(x_i)$ can be arranged as $( x_{m\alpha i})$ with
\begin{equation*}
 m=1,2,\dots, \quad \alpha=1,
\dots, N_m, \quad i=1,\dots,m
\end{equation*}
such that\hfill\break
(a) $d_{m \alpha i}
=m+1-2i$ for $d_{m\alpha i} := \ord(x_{m\alpha i})$.
\hfill\break
(b)
$z_{m\alpha1}^{0}=z_{m\alpha2}^{0}=
\cdots = z_{m\alpha m}^{0}$,
where $z_{m\alpha i}=\tilde{x}_{m\alpha i}$.
\end{definition}
For each $1 \le \alpha \le N_m$,
$(z_{m\alpha i})_{i=1}^m$  is called an $m$-string.
$N_m$ is the number of $m$-strings.
When considering string solutions, we
denote the quantity in (b) by $z^0_{m \alpha}$, and call it the
{\em string center}.
We set
\begin{equation*}
q^{\zeta_{m\alpha i}}y_{m\alpha i}(q) = z_{m\alpha i}(q) - z_{m\alpha i-1}(q)
\qquad 2 \le i \le m,
\end{equation*}
where $\zeta_{m\alpha i} = \ord(z_{m\alpha i} - z_{m\alpha i-1})
\in \Z_{\ge 1}$.
For a string solution of pattern $N$,
the Bethe equation (\ref{eq:ba0}) reads
\begin{equation}\label{eq:ba2}
F_{m\alpha i+}G_{m\alpha i-}=-
F_{m \alpha i-}G_{m \alpha i+},
\end{equation}
where
\begin{alignat*}{2}
F_{m \alpha i+}&=
\prod_{s\geq 1}(z_{m \alpha i}q^{d_{m \alpha i}+{s}}
-1)^{\nu_s},
&\ 
G_{m \alpha i+}&=
\prod_{k\beta j}
(z_{m \alpha i}q^{d_{m \alpha i}+2}
 -z_{k \beta j}q^{d_{k\beta j}}),
\\
F_{m \alpha i-}&=
\prod_{s\geq 1}(z_{m \alpha i}q^{d_{m \alpha i}}
-q^{s})^{\nu_s},
&\
G_{m\alpha i-}&=
\prod_{k\beta j}
(z_{m \alpha i}q^{d_{m \alpha i}}
 -z_{k \beta j}q^{d_{k \beta j}+2}).
\end{alignat*}
Here $\prod_{k\beta j}$ means
$\prod_{k \ge 1} \prod_{\beta = 1}^{N_k} \prod_{j = 1}^k$.
Let us  extract the factors $y_{m\alpha i}$ from
$G_{m\alpha i\pm}$
by introducing $G^{\prime}_{m\alpha i\pm}$ as follows:
\begin{align*}
G_{m\alpha i+} &=
\begin{cases}
G^{\prime}_{m\alpha 1+} & i=1\\
G^{\prime}_{m\alpha i+}q^{d_{m\alpha i}+2+\zeta_{m\alpha i}}y_{m\alpha i}
&2 \le i \le m,
\end{cases}\\
G_{m\alpha i-} &=
\begin{cases}
G^{\prime}_{m\alpha i-}(-q^{d_{m\alpha i}+\zeta_{m\alpha i+1}}y_{m\alpha i+1})
&1 \le i \le m-1\\
-G^{\prime}_{m\alpha m-} & i=m.
\end{cases}
\end{align*}
Now the Bethe equation (\ref{eq:ba2}) takes the form:
\begin{alignat}{3}
\tilde{F}_{m\alpha 1+}
\tilde{G}^{\prime}_{m\alpha 1-}y_{m\alpha 2}
& =
\tilde{F}_{m\alpha 1-}\tilde{G}^{\prime}_{m\alpha 1+}
&\qquad &i =
1,\label{eq:fgy11}\\
\tilde{F}_{m\alpha i+}
\tilde{G}^{\prime}_{m\alpha i-}y_{m\alpha i+1}
& =
\tilde{F}_{m\alpha i-}\tilde{G}^{\prime}_{m\alpha i+}y_{m\alpha i}
&\qquad &2 \le i \le m-1,\label{eq:fgy22}\\
\tilde{F}_{m\alpha m+}\tilde{G}^{\prime}_{m\alpha m-}
& =
\tilde{F}_{m\alpha m-}
\tilde{G}^{\prime}_{m\alpha m+}y_{m\alpha m}
&\qquad &i = m.\label{eq:fgy33}
\end{alignat}
In particular,
\begin{equation}\label{eq:ratio2}
1 = (-1)^m \prod_{i=1}^m
\frac{F_{m\alpha i+}G_{m\alpha i-}}{F_{m\alpha i-}G_{m\alpha i+}} =
\prod_{i=1}^m\frac{\tilde{F}_{m\alpha i+}\tilde{G}^{\prime}_{m\alpha i-}}
{\tilde{F}_{m\alpha i-}\tilde{G}^{\prime}_{m\alpha i+}},
\end{equation}
where the latter equality is the identity
as  holomorphic functions around $q=0$.

\subsection{$q\rightarrow 0$ limit of Bethe equation}
\label{subsec:q0be}

Suppose that $(x_{m\alpha i})$ is a string solution to the Bethe equation
(\ref{eq:ba2}).
Taking the leading coefficients of (\ref{eq:fgy11})--(\ref{eq:fgy33})
we get the $q \rightarrow 0$ limit:
\begin{alignat}{3}
F^0_{m\alpha 1+}G^{\prime \; 0}_{m\alpha 1-}y^0_{m\alpha 2}
& =
F^0_{m\alpha 1-}G^{\prime \; 0}_{m\alpha 1+}
&\qquad &i = 1,\label{eq:fgy1}\\
F^0_{m\alpha i+}G^{\prime \; 0}_{m\alpha i-}y^0_{m\alpha i+1}& =
F^0_{m\alpha i-}G^{\prime \; 0}_{m\alpha i+}y^0_{m\alpha i}
&\qquad &2 \le i \le m-1,\label{eq:fgy2}\\
F^0_{m\alpha m+}G^{\prime \; 0}_{m\alpha m-} &=
F^0_{m\alpha m-}G^{\prime \; 0}_{m\alpha m+}y^0_{m\alpha m}
&\qquad &i = m.\label{eq:fgy3}
\end{alignat}
In particular,
\begin{equation}\label{eq:ratio1}
1 = (-1)^m \prod_{i=1}^m
\frac{F^0_{m\alpha i+}G^0_{m\alpha i-}}{F^0_{m\alpha i-}G^0_{m\alpha i+}} =
\prod_{i=1}^m\frac{F^0_{m\alpha i+}G^{\prime \; 0}_{m\alpha i-}}
{F^0_{m\alpha i-}G^{\prime \; 0}_{m\alpha i+}}
\end{equation}
holds for the leading coefficients.

\subsection{Generic string solution}

In order to estimate the order of the Bethe equation (\ref{eq:ba2})
we introduce
\begin{displaymath}
\begin{split}
\xi_{m\alpha i+}&=
\sum_{s\geq 1}
\nu_s \min(m+1-2i+s,0),\\
\xi_{m\alpha i-}&=
\sum_{s\geq 1}
\nu_s \min(m+1-2i,s),\\
\eta_{m\alpha i+}&=
\sum_{k\beta j}
\min(m+3-2i,k+1-2j),\\
\eta_{m\alpha i-}&=
\sum_{k\beta j}
\min(m+1-2i,k+3-2j),
\end{split}
\end{displaymath}
where $\sum_{k\beta j}$ is the abbreviation of
$\sum_{k \ge 1} \sum_{\beta = 1}^{N_k} \sum_{j = 1}^k$.
In general one has  $\xi_{m\alpha i\pm} \le \ord(F_{m\alpha i\pm})$,
$\eta_{m\alpha i+} \le \ord(G_{m\alpha i+})+(1-\delta_{i,1})\zeta_{m\alpha i}$
and 
$\eta_{m\alpha i-} \le \ord(G_{m\alpha i-})+(1-\delta_{i,m})\zeta_{m\alpha
i+1}$.
Let us consider the simplest situation when these inequalities are saturated.
\begin{definition}\label{def:generic1}\upshape
A string solution $(x_{m\alpha i})$ to (\ref{eq:ba0})
is called {\em generic\/} if the following
condition is valid:
\begin{align*}
\ord(F_{m\alpha i\pm}) &= \xi_{m\alpha i\pm},\\
\ord(G_{m\alpha i+}) &=
\begin{cases}
\eta_{m\alpha 1+} & i = 1\\
\eta_{m\alpha i+} + \zeta_{m\alpha i} & 2 \le i \le m,
\end{cases}\\
\ord(G_{m\alpha i-})&=
\begin{cases}
\eta_{m\alpha i-} + \zeta_{m\alpha i+1} & 1 \le i \le m-1\\
\eta_{m\alpha m-} & i = m.
\end{cases}
\end{align*}
\end{definition}

Our results in the rest of Section \ref{sec:beq0} mostly concern
generic string solutions.
Given the quantum space data $\nu =(\nu_s)$ and the string pattern
$N = (N_m)$, we put
\begin{align}
P_m &= P_m(\nu,N) = \gamma_m - 2\sum_{k \ge 1} \min(m,k)N_k,
\label{eq:pdef}\\
\gamma_m &= \gamma_m(\nu) = \sum_{k\ge 1} \min(m,k) \nu_k.\label{eq:gammadef}
\end{align}
\begin{lemma}\label{lem:order1}
For $1 \le i \le m$ we have
\begin{align*}
& (\xi_{m\alpha i+}+\eta_{m\alpha i-})
-(\xi_{m\alpha i-}+\eta_{m\alpha i+})\\
&\qquad =
\begin{cases}
-(P_{m+1-2i}+N_{m+1-2i}) &
i < {m+1\over 2}\\
0 & i = {m+1\over 2}\\
P_{2i-m-1}+N_{2i-m-1} &
i > {m+1\over 2}.
\end{cases}
\end{align*}
\end{lemma}
\begin{proposition}\label{pr:order}
A necessary condition for the existence of
a ge\-ner\-ic
 string solution of pattern $N = (N_m)$ is as follows:
If $N_m \ge 1$, then
\begin{displaymath}
\begin{split}
&P_{m-1}+N_{m-1} \ge 1,\\
&(P_{m-1}+N_{m-1}) + (P_{m-3}+N_{m-3})   \ge 1,\\
&\qquad \cdots \\
&(P_{m-1}+N_{m-1}) + (P_{m-3}+N_{m-3}) + \cdots +
\begin{cases}
(P_1 + N_1) \ge 1 & m:\text{even}\\
(P_2 + N_2) \ge 1 & m:\text{odd}.
\end{cases}
\end{split}
\end{displaymath}
\end{proposition}
\begin{proof}
{}From the condition 
$\ord(F_{m\alpha i+}G_{m\alpha i-}) = \ord(F_{m\alpha i-}G_{m\alpha i+})$,
we have
\begin{displaymath}
 (\xi_{m\alpha i+}+\eta_{m\alpha i-})
-(\xi_{m\alpha i-}+\eta_{m\alpha i+})
=
\begin{cases}
-\zeta_{m\alpha 2} & i = 1\\
\zeta_{m\alpha i} - \zeta_{m\alpha i+1}& 2 \le i \le m-1\\
\zeta_{m\alpha m}& i = m.
\end{cases}
\end{displaymath}
Solving this by using {\sc Lemma} \ref{lem:order1}, we get
\begin{displaymath}
\zeta_{m\alpha i} = \zeta_{m\alpha \,m+2-i} =
\sum_{k=1}^{\min(i-1,m+1-i)}(P_{m+1-2k} +
N_{m+1-2k})
\qquad 2 \le i \le m.
\end{displaymath}
In order to have
$\zeta_{m\alpha 2}$, $\zeta_{m\alpha 3}$, \dots, $\zeta_{m\alpha m}
\ge 1$, the condition in the proposition must hold.
\end{proof}
\par
In the remainder of this section we shall exclusively consider the case
$\forall P_m(\nu,N) \ge 0$.
Thus the necessary condition in {\sc Proposition} \ref{pr:order}
reduces to
\begin{equation}\label{eq:ncondition}
\text{$P_m \ge 0$ \ for any $m$,\quad
 $P_{m-1} + N_{m-1} \ge 1$
 \ whenever \ $N_m \ge 1$.}
\end{equation}

\subsection{String center equation (SCE)}\label{subsec:sce}

\begin{theorem}\label{th:A}
Let $(x_{m\alpha i})$ be a generic string solution of pattern $N$.
Then its string center $(z^0_{m\alpha})$ satisfies
the equation:
\begin{gather}
\prod_{k\ge 1} \prod_{\beta = 1}^{N_k}
(z^0_{k\beta})^{A_{m\alpha,k\beta}} = (-1)^{P_m+N_m+1}
\qquad m \ge 1, \; 1 \le \alpha \le N_m,\label{eq:sce1}\\
A_{m\alpha, k\beta} := \delta_{m k}\delta_{\alpha \beta}(P_m+N_m) +
2\min(m,k) - \delta_{m k}.\label{eq:amat}
\end{gather}
\end{theorem}
We call (\ref{eq:sce1}) the string center equation (SCE).
It is a linear congruence equation in the sense of
(\ref{eq:sce2}).
\begin{proof}
Let us compute the ratio (\ref{eq:ratio1}) explicitly.
\begin{gather*}
\prod_{i=1}^m F_{m\alpha i\epsilon }^0 =
\begin{cases}
(-1)^{\gamma_m}
\prod_{s\ge 1}(f^{(s)}_{m\alpha})^{\nu_s}& \epsilon = +\\
(z^0_{m\alpha})^{\gamma_m}
\prod_{s\ge 1}(f^{(s)}_{m\alpha})^{\nu_s}& \epsilon = -,
\end{cases}\\
f^{(s)}_{m\alpha} = \begin{cases}
1 & m \le s\\
(-z^0_{m\alpha})^{\frac{m-s}{2}} &
m > s, s \equiv m\ \mathrm{mod}\ 2\\
(-z^0_{m\alpha})^{\frac{m-s-1}{2}}(z^0_{m\alpha}-1) &
m > s, s \not \equiv m\ \mathrm{mod}\ 2.
\end{cases}
\end{gather*}
In  order to  calculate $\prod_{i=1}^m(G^0_{m\alpha i-}/G^0_{m\alpha i+})$,
it is convenient first to evaluate
\begin{gather*}
\begin{split}
&\prod_{i=1}^m\prod_{j=1}^k(z_{m\alpha i}q^{d_{m\alpha i}+1
+\epsilon} -z_{k\beta j}q^{d_{k\beta j}+1-\epsilon})^0\\
&\qquad =
\begin{cases}
(-z^0_{k\beta})^{2\min(m,k)-\delta_{m k}}g^{k\beta}_{m\alpha} &
\epsilon = 1\\
(z^0_{m\alpha})^{2\min(m,k)-\delta_{m k}}
(-1)^{(m-1)\delta_{m k}\delta_{\alpha \beta}}g^{k\beta}_{m\alpha} 
& \epsilon = -1,
\end{cases}\\
\end{split}\\
g^{k\beta}_{m\alpha} =
\begin{cases}
(-z^0_{m\alpha}z^0_{k\beta})^{\frac{mk}{2}-\min(m,k)}
& m \not\equiv k\ \mathrm{mod}\ 2\\
(-z^0_{m\alpha}z^0_{k\beta})^{\frac{mk}{2}-\frac{3}{2}
\min(m,k)+\delta_{m k}}
& m\alpha \neq k\beta,
 \ m \equiv k \ \mathrm{mod}\ 2\\
\quad
\times (z^0_{m\alpha}-z^0_{k\beta})^{\min(m,k)-\delta_{m k}}
& \\
(-z^0_{m\alpha}z^0_{k\beta})^{\frac{mk}{2}-
\frac{3}{2}\min(m,k)+ \delta_{m k}}
& m\alpha = k\beta.\\
\quad \times y^0_{m\alpha 2} \cdots y^0_{m\alpha m}
&
\end{cases}
\end{gather*}
The factors $\prod_{s\ge 1}(f^{(s)}_{m\alpha})^{\nu_s}$ and
$g^{k\beta}_{m\alpha}$ are
all nonzero for a generic string solution $(x_{m\alpha i})$.
They are canceled in the ratio (\ref{eq:ratio1}) and we find
\begin{equation}
\label{eq:ratio}
1 = (-1)^m \prod_{i=1}^m
\frac{F^0_{m\alpha i+}G^0_{m\alpha i-}}{F^0_{m\alpha i-}G^0_{m\alpha i+}}
= (-1)^{P_m+N_m+1} \prod_{k\beta}(z^0_{k\beta})^{-A_{m\alpha,k\beta}}.
\end{equation}
\end{proof}

\begin{remark}\label{rem:generic1}\upshape
{}From the condition 
$\prod_{s\ge 1}(f^{(s)}_{m\alpha})^{\nu_s}$, $g^{k\beta}_{m\alpha}\neq 0$
in the above proof,
we see that a string solution $(x_{m\alpha i})$ is  generic if and only if
\begin{equation}\label{eq:nagoya}
\begin{gathered}
\prod_{1 \le s < m,\, s \not\equiv m (2)}(z^0_{m\alpha} - 1 )^{\nu_s} \neq 0,\\
\prod_{k \beta  (\neq m\alpha),\, k \equiv m (2)}
(z^0_{m\alpha} - z^0_{k\beta})^{\min(m,k)-\delta_{m k}}
 \neq 0\;
\end{gathered}
\end{equation}
for any $m \ge 1,\; 1 \le \alpha \le N_m$.
In the latter, the power $\min(m,k)-\delta_{m k}$ implies that
the collision $z^0_{m\alpha} = z^0_{k\beta}$ is allowed only when $m = k = 1$.
\end{remark}

\begin{definition}\label{def:generic2}\upshape
A solution to SCE (\ref{eq:sce1}) is called {\em generic\/}
if it satisfies the  condition (\ref{eq:nagoya}).
\end{definition}
By the definition the solutions to SCE (\ref{eq:sce1}) arising
{}from the string centers of
generic string solutions as in {\sc Theorem} \ref{th:A} are generic.

SCE (\ref{eq:sce1})
becomes the linear congruence equation
in terms of the variables $u_{k\beta} \in \R/\Z$ defined by
$z^0_{k\beta} = \exp(2\pi\sqrt{-1}u_{k\beta})$:
\begin{equation}\label{eq:sce2}
\sum_{k\ge 1}\sum_{\beta=1}^{N_k} A_{m\alpha,k\beta} u_{k\beta} \equiv
\frac{P_m + N_m + 1}{2} \quad \mathrm{mod}\ \Z.
\end{equation}
This will also be called SCE.
In the limit $\hbar \rightarrow \infty$, 
the asymptotic behavior of the original variable
$u_j = \frac{1}{2\pi\sqrt{-1}}\log(x_{k\beta i})$ 
in (\ref{eq:be}) is $u_{k\beta}+ \sqrt{-1}\hbar(k+1-2i)$.

\subsection{Lifting generic solutions of SCE
to generic string solutions}\label{subsec:lift}

In Section \ref{subsec:sce} we have seen that string
centers of a generic string solution
to the Bethe equation (\ref{eq:ba2}) yield a generic solution to SCE.
Here we show the inverse.
\begin{theorem}\label{th:B}
Suppose that $N = (N_m)$ obeys the condition (\ref{eq:ncondition}).
Let $(z^0_{m\alpha})$ be a generic solution to SCE (\ref{eq:sce1}).
Then there exists a unique generic string solution $(x_{m\alpha i}(q))$ to the
Bethe equation (\ref{eq:ba0}) such that
$z_{m\alpha i}(0) = z^0_{m\alpha}\; (1 \le i \le m)$.
\end{theorem}
Define the variables (holomorphic functions of $q$) $w_{m\alpha i}$ by
\begin{equation*}
w_{m\alpha i} =
\begin{cases}
z_{m\alpha i} & i=1\\
y_{m\alpha i} & 2 \le i \le m.
\end{cases}
\end{equation*}
Then 
\begin{equation}\label{eq:wxyz}
z_{m\alpha i} = w_{m\alpha 1} + q^{\zeta_{m\alpha 2}}w_{m\alpha 2} + \cdots
+ q^{\zeta_{m\alpha i}}w_{m\alpha i}\qquad 1 \le i \le m.
\end{equation}
Denote the $i$th equation of (\ref{eq:fgy11})--(\ref{eq:fgy33})
by $L_{m\alpha i} = R_{m\alpha i}$.
Let $A = (A_{m\alpha, k\beta})$ be the matrix of size
$N_1 + N_2 + \cdots$ defined by (\ref{eq:amat}).
Similarly let $J = (J_{m\alpha i, k\beta j})$ be the matrix of size
$N_1 + 2N_2 + \cdots (=M)$ defined by
$J_{m\alpha i, k\beta j} = \frac{\partial}{\partial w_{k\beta j}}
\left(\frac{L_{m\alpha i}}{R_{m\alpha i}}-1\right)$.
\begin{lemma}\label{lem:jacobian}
Suppose (\ref{eq:ncondition}) is satisfied.
Then  $\det J $ is nonzero at $q=0$
(i.e., $\det(J^0_{m\alpha i, k\beta j}) \neq 0$)
if $\det A \neq 0$.
\end{lemma}
\begin{proof}
Owing to the assumption (\ref{eq:ncondition}), we have
$\forall \zeta_{m\alpha i}  \ge 1$.
Since $L^0_{m\alpha i} = R^0_{m\alpha i} \neq 0$,
it suffices to show $\det{\mathcal{J}} \neq 0$ for 
${\mathcal{J}}_{m\alpha i, k\beta j} =
\frac{\partial}{\partial w_{k\beta j}}\log
\frac{L_{m\alpha i}}{R_{m\alpha i}}
$.
{}From (\ref{eq:wxyz}) both
$\frac{\partial\tilde{F}_{m\alpha i \pm}}
{\partial w_{k\beta j}}$ and
$\frac{\partial\tilde{G}^{\prime}_{m\alpha i \pm}}
{\partial w_{k\beta j}}$
for  $j \neq 1$ are zero at $q=0$.
Thus among ${\mathcal{J}}^0_{m\alpha i, k\beta j}$'s the non-vanishing
 ones are only
${\mathcal{J}}^0_{m\alpha i, k\beta 1}$ $(1\leq i \leq m)$,
${\mathcal{J}}^0_{m\alpha i, m\alpha i}
 = -1/{y^0_{m\alpha i}}$ $(2 \le i \le m)$,
and 
${\mathcal{J}}^0_{m\alpha i, m\alpha i+1}
= 1/{y^0_{m\alpha i+1}}$ $(1 \le i \le m-1)$.
%
Let $\vec{{\mathcal{J}}}^{\;0}_{m\alpha i} = ({\mathcal{J}}^0_{m\alpha i,
k\beta j})_{k\beta j}$
be a row vector of the matrix ${\mathcal{J}}$.
In view of the above result,
the linear dependence $\sum_{m\alpha i}c_{m\alpha
i}\vec{{\mathcal{J}}}^{\;0}_{m\alpha i} = 0$
can possibly holds only when $c_{m\alpha i}$ is independent of $i$.
Consequently we consider the equation
$\sum_{m\alpha}c_{m\alpha}\sum_{i=1}^m \vec{{\mathcal{J}}}^{\;0}_{m\alpha
i} = 0$.
The $(k\beta 1)$-th component of the vector
$\sum_{i=1}^m \vec{{\mathcal{J}}}^{\;0}_{m\alpha i}$
is given by
\begin{displaymath}
\lim_{q \rightarrow 0} \frac{\partial}{\partial w_{k\beta 1}} \log
\prod_{i=1}^m \frac{\tilde{F}_{m\alpha i+}\tilde{G}^{\prime}_{m\alpha i-}}
{\tilde{F}_{m\alpha i-}\tilde{G}^{\prime}_{m\alpha i+}}
= \frac{\partial}{\partial z^0_{k\beta}} \log
\prod_{i=1}^m \frac{-F^0_{m\alpha i+}G^0_{m\alpha i-}}
{F^0_{m\alpha i-}G^0_{m\alpha i+}},
\end{displaymath}
where we have taken (\ref{eq:ratio2}),  (\ref{eq:wxyz}) and
$\forall \zeta_{m\alpha i} \ge 1$ into account.
Due to (\ref{eq:ratio}) the last expression is equal to
$-A_{m\alpha, k\beta}/z^0_{k\beta}$.
Therefore the equation
$\sum_{m\alpha}c_{m\alpha} \sum_{i=1}^m \vec{{\mathcal{J}}}^{\;0}_{m\alpha
i} = 0$
is equivalent to
$\sum_{m\alpha}c_{m\alpha}A_{m\alpha, k\beta} = 0$ for any $k \beta$.
This admits only the trivial solution for
$c_{m\alpha}$ if $\det A \neq 0$.
\end{proof}

{\sc Proof of Theorem \ref{th:B}}.
The Bethe equations (\ref{eq:fgy11})--(\ref{eq:fgy33})
are simultaneous equations on the variables
$(w_{m\alpha i}, q)$.
At $q=0$, (\ref{eq:fgy11})--(\ref{eq:fgy33})
reduce to (\ref{eq:fgy1})--(\ref{eq:fgy3}).
The latter fix $(y^0_{m\alpha i})$ unambiguously
once a generic solution $(z^0_{m\alpha})$ to SCE is given.
Denote the resulting value of $w_{m\alpha i}$ by
$w^0_{m\alpha i}$.
Thus (\ref{eq:fgy11})--(\ref{eq:fgy33}) are valid
at $(w_{m\alpha i},q) = (w^0_{m\alpha i},0)$.
{}From the implicit function theorem, there uniquely
exist the functions $w_{m\alpha i}(q)$ satisfying
(\ref{eq:fgy11})--(\ref{eq:fgy33}) and
$w_{m\alpha i}(0) = w^0_{m\alpha i}$,
if the Jacobian
$\det J$ at $(w_{m\alpha i},q) = (w^0_{m\alpha i},0)$
is nonzero.
By {\sc Lemma} \ref{lem:jacobian} and
{\sc Corollary} \ref{cor:detanonzero}
this has been guaranteed under
(\ref{eq:ncondition}).
\begin{flushright}
$\square$
\end{flushright}

{}From {\sc Theorem} \ref{th:A} and {\sc Theorem} \ref{th:B}
we have
\begin{corollary}\label{cor:AB}
Suppose the pattern $N = (N_m)$ satisfies
 (\ref{eq:ncondition}).
Then there is a one-to-one correspondence between
generic string solutions to the Bethe equation (\ref{eq:ba0})
and generic solutions to SCE (\ref{eq:sce1}).
\end{corollary}

\subsection{Off-diagonal solution}\label{subsec:off}

In {\sc Corollary} \ref{cor:AB} 
to restrict oneself to the admissible solutions to the Bethe equation 
is natural, because
otherwise the associated Bethe vectors are vanishing \cite{TV}.
On the other hand, 
the limitation to the generic case has been made by a technical reason.
In fact the most essential constraint on the solution of the Bethe
equation (\ref{eq:ba0}) is that $x_1$, \dots, $x_M$ are all
distinct.
Otherwise the Bethe vectors are again vanishing.
To make it precise in the present context, we introduce
\begin{definition}\label{def:off1}\upshape
A meromorphic solution $(x_i)$
of (\ref{eq:ba0})
is called {\em diagonal (off-diagonal)\/} if
$x_i=x_j$ for some $i\neq j$
 as a function of $q$  around $q=0$ (otherwise).
\end{definition}
\begin{definition}\label{def:off2}\upshape
A solution $(z_{m\alpha}^0)$ of SCE
is called {\em diagonal (off-diagonal)\/} if
$z^0_{m\alpha} = z^0_{m\beta}$ for some $\alpha\neq\beta$
(otherwise).
\end{definition}
According to the above definition,
the string solution $(z_{m\alpha i})$ of the Bethe
equation is {\em diagonal (off-diagonal)\/} if
$z_{m\alpha i}(q)= z_{k\beta j}(q)$ and
$d_{m\alpha i}=d_{k \beta j}$ for some
$m\alpha i \neq k\beta j$ (otherwise).
We conjecture that the number of off-diagonal string solutions of the Bethe
equation is equal to the number of off-diagonal solutions of SCE
under a certain condition like $(\ref{eq:ncondition})$.
So far we have been unable to solve the discrepancy between ``generic" and
``off-diagonal".
\begin{remark}
In \cite{LS} SCE has been
given for the XXZ case ($\nu_s = L\delta_{1 s}$).
It is claimed that SCE is satisfied irrespective of whether
a string solution is generic or not.
Moreover, each off-diagonal solution of SCE gives rise to
an off-diagonal string solution.
\end{remark}


\section{Counting of off-diagonal solutions}\label{sec:counting}
This section is devoted to an expository proof of
{\sc Theorem} \ref{th:Rnds}.
It provides an explicit combinatorial formula
counting the number of off-diagonal solutions to
SCE (\ref{eq:sce1}) in the sense of Section \ref{subsec:off}.
We assume that the quantum space data
$\nu =(\nu_s)$ and the string pattern $N = (N_m)$ satisfy the
condition
\begin{equation}\label{eq:ppositive}
\text{$P_m(\nu,N) \ge 0$ \ whenever \ $N_m \ge 1$}
\end{equation}
throughout Section \ref{sec:counting}.
This is a milder condition than $\forall P_m \ge 0$.

\subsection{Rule of counting} \label{subsec:rule}

We shall work with the logarithmic form of SCE (\ref{eq:sce2}) presented as
\begin{equation}\label{eq:sce3}
A \vec{u} \equiv \vec{c} \quad \mathrm{mod}\ \Z^d,
\end{equation}
where $d = N_1 + N_2 + \cdots$.
Here $\vec{u} = (u_{k\beta})$ is the unknown and $\vec{c}$
is some constant vector.
The $d$-dimensional matrix $A = (A_{m\alpha, k\beta})$ is
specified by (\ref{eq:amat}).
It has a block structure according to the string pattern.
For example,  if only $N_1, N_2$ and $N_3$ are non-zero,
it looks as (${\mathcal P}_i = P_i + N_i$)
{\small
\begin{displaymath}
\left(
\begin{array}{cccccccccccc}
{\mathcal P}_1+1 &  \cdots & 1 &
2 & \cdots & 2 &
2 & \cdots & 2 \\
\vdots &\ddots  & \vdots &
\vdots &   & \vdots  &
\vdots &   & \vdots  \\
1 & \cdots  & {\mathcal P}_1+1 &
2 & \cdots & 2 &
2 & \cdots & 2 \\
2 & \cdots & 2 &
{\mathcal P}_2+3 & \cdots & 3 &
4 & \cdots & 4 \\
\vdots &   & \vdots  &
\vdots & \ddots & \vdots &
\vdots &   & \vdots  \\
2 & \cdots & 2 &
3 & \cdots & {\mathcal P}_2+3 &
4 & \cdots  & 4 \\
2 & \cdots & 2 &
4 & \cdots & 4 &
{\mathcal P}_3+5 &  \cdots & 5 \\
\vdots &   & \vdots  &
\vdots & \ddots  & \vdots  &
\vdots &   & \vdots \\
2 & \cdots  & 2 &
4 & \cdots  & 4 &
5 & \cdots & {\mathcal P}_3+5
\end{array}\right),
\end{displaymath}
}
which consists of 9 sub-matrices of size $N_i \times N_j$
($1 \le i,j \le 3$).

We will seek the number of off-diagonal solutions to SCE (\ref{eq:sce3})
in the sense of Section \ref{subsec:off}.
Let us remember the three essential rules to count it.
First, there is no change in (\ref{eq:be}) under any
integer shift of $u_j$.
Accordingly we should not distinguish the solutions
$\vec{u}$ and $\vec{u}^{\, \prime}$ to (\ref{eq:sce3})
if $\vec{u} - \vec{u}^{\, \prime} \in \Z^d$.
In other words, we consider $u_{k \beta} \in \R/\Z$ rather than $\R$.
Secondly,  the original Bethe equation (\ref{eq:be}) is symmetric with
respect to
$u_1$, \dots, $u_M$, but their permutation does not lead to a new
Bethe vector as is well known.
Consequently, for a solution  $\vec{u} = (u_{k\beta})$
of  pattern  $N = (N_m)$,
we should regard
\begin{equation*}
(u_{k 1}, u_{k2}, \ldots, u_{kN_k}) \in
\left(\R/\Z\right)^{N_k} /{\frak S}_{N_k}
\end{equation*}
for each $k \in \N$.
Here ${\frak S}_{N_k}$ stands for the symmetric group
permuting the  $N_k$ components.
Last but most importantly, {\sc Definition} \ref{def:off2}
postulates that
$u_{k\beta} \neq u_{k\beta'} \in \R/\Z$ if
$1 \le \beta \neq \beta' \le N_k$ for each $k \in \N$.
To summarize these rules, we start with a fixed string pattern $(N_m)$ and
specify $A$ and $\vec{c}$ by  (\ref{eq:sce2}).
Then we are to count the number of solutions
$\vec{u} = (u_{k \beta})$ to
SCE  (\ref{eq:sce3}) such that
\begin{gather*}
(u_{k1}, u_{k2}, \ldots, u_{kN_k}) \in
\bigl(\left(\R/\Z\right)^{N_k} - \Delta_{N_k} \bigr) /{\frak S}_{N_k},\\
\Delta_n  = \{(v_1, \ldots, v_n)
\in (\R/\Z)^n \mid v_\alpha = v_\beta\ \text{for some}\
 1 \le \alpha \neq
\beta \le n\}
\end{gather*}
for each $k \in \N$.
In practice, one just has to find the number of solutions
such that $(u_{k1}, u_{k2}, \ldots, u_{kN_k}) \in
\left(\R/\Z\right)^{N_k} - \Delta_{N_k}$, and
divide afterwards by $N_k!$ for each $k$.

\subsection{Example}\label{subsec:examples}
Before treating the general case,
let us illustrate how to enumerate the off-diagonal solutions with
a simplest example $\nu_s = L\delta_{s,1}$.
This corresponds to the spin $\frac{1}{2}$ XXZ model with $L$-sites.
Thus the character of the quantum space is expanded
 as ($x = e^{\Lambda_1}$)
\begin{equation*}
(x+x^{-1})^L = x^L + \binom{L}{1}x^{L-2} +
\binom{L}{2}x^{L-4} + \binom{L}{3}x^{L-6} + \cdots,
\end{equation*}
according to the `magnon number' $M = 0$, $1$, $2$, $3$,
 \dots.
Here and in what follows the symbol $\binom{\,\cdot\,}{\,\cdot\,}$ will always
denote the generalized binomial coefficient:
\begin{equation*}
\binom{\xi}{n} = \begin{cases}
\xi(\xi-1) \cdots (\xi-n+1)/n! & \text{ if } n \in \Z_{\ge 1}\\
1 & \text{ if  } n = 0\\
0 & \text{ otherwise},
\end{cases} \qquad \xi \in \C.
\end{equation*}
The counting of the off-diagonal solutions
sketched below up to $M = 3$ is useful to gain
the idea for treating the general case.
The result
coincides with the weight multiplicity $\binom{L}{M}$
appearing in the above expansion.
The proof of the coincidence in the general case will be given in
Section \ref{sec:multiplicity}.
According  to (\ref{eq:ppositive})  we assume that
$P_m = L - 2\sum_{k\ge 1}\min(m,k)N_k \ge 0$ whenever $N_m > 0$.
We use the following fact:
\begin{proposition}\label{lem:basic}
Let $B$ be an $n$ by $n$ integer matrix with $\det B \neq 0$.
Then for any $\vec{b} \in \R^n$, the equation
$B \vec{x} \equiv \vec{b}$ $\mathrm{mod}$ $\Z^n$ has exactly 
$\vert \det B \vert$ solutions $\vec{x}$ in $(\R/\Z)^n$.
\end{proposition}
\begin{proof}
Because of linearity,
 it is enough to show the statement for $\vec{b}=\vec{0}$.
Consider a pair of lattices
$L \supset M$, where $L$ consists of the solutions
$\vec{x}\in \mathbb{R}^n$ for the
homogeneous
equation $B\vec{x}\equiv\vec{0}$ mod $\mathbb{Z}^n$,
and $M=\mathbb{Z}^n$. The column vectors
 $\vec{f}_1,\dots,\vec{f}_n$ of
$B^{-1}$ is a basis of $L$.
It is well-known \cite[1.2.2, Lemma 1]{C}
that 
$|L/M|$ equals to $|\det(\vec{f}_1,\dots,\vec{f}_n)|^{-1}$,
which is $|\det B|$.
\end{proof}

%
As it turns out, the condition (\ref{eq:ppositive})
assures $\det B > 0$ not only for $B = A$ but also
for all the relevant matrices $B$ coming into the game.
(cf.\ {\sc Lemma} \ref{lem:positive}.)

{\em $M = 0$ Case}.
The equation (\ref{eq:sce3}) is void.
We just define the number of off-diagonal solutions  to be 1
in agreement with $\binom{L}{0}$.

{\em $M = 1$ Case}.
The only string pattern (\ref{eq:mpartition})
is $N_1= 1$.
Then (\ref{eq:sce3}) is a scalar equation $L u_{1,1} \equiv c$ mod $\Z$
for some $c \in \R$.
Thus the number of the off-diagonal solution is $L$
in agreement with $\binom{L}{1}$.

{\em $M = 2$ Case}.
There are two string patterns, (i) $N_2=1$ and (ii) $N_1 = 2$,
that satisfy (\ref{eq:mpartition}).
(i) The equation (\ref{eq:sce3}) is
$L u_{2,1} \equiv c$ for some  $c \in \R$.
Thus there are $L$ off-diagonal solutions.
(ii) The equation (\ref{eq:sce3}) reads
\begin{equation*}
A \vec{u} =
\begin{pmatrix}
L-1 & 1 \\
1 & L-1
\end{pmatrix}
\begin{pmatrix} u_{1,1} \\ u_{1,2} \end{pmatrix}
\equiv
\vec{c}\quad \mathrm{mod}\ \Z^2
\end{equation*}
for some $\vec{c}$.
{}From the assumption $P_1 = L - 4 \ge 0$,
there are $\det A = L(L-2)$  solutions to this by 
{\sc Proposition} \ref{lem:basic}.
However they contain the diagonal ones $u_{1,1} = u_{1,2}$.
Under this constraint the above matrix equation reduces to
$L u_{1,1} \equiv c'$ with some $c'$, telling that the number of
the diagonal solutions is $L$.
Therefore the non-diagonal solutions from (ii) is enumerated as
$(L(L-2) - L)/2$ by recalling the $\frak{S}_2$ redundancy.
Collecting the contributions from the string patterns (i) and (ii) one finds
the number of off-diagonal solutions
\begin{equation*}
L + \frac{1}{2}(L(L-2) - L) = \binom{L}{2}.
\end{equation*}

{\em $M = 3$ Case}.
There are three string patterns, (i) $N_3=1$,
 (ii) $N_1 = N_2 = 1$
and (iii) $N_1 = 3$
that satisfy (\ref{eq:mpartition}).
(i) The equation (\ref{eq:sce3}) is
$L u_{3,1} \equiv c$ for some  $c \in \R$.
Thus there are $L$ off-diagonal solutions.
(ii) The equation (\ref{eq:sce3}) reads
\begin{equation*}
A \vec{u} =
\begin{pmatrix}
L-2 & 2 \\
2 & L-2
\end{pmatrix}
\begin{pmatrix} u_{1,1} \\ u_{2,1} \end{pmatrix}
\equiv
\vec{c}\quad \mathrm{mod}\ \Z^2
\end{equation*}
for some $\vec{c}$.
{}From the assumption $P_1 = L - 4, P_2 = L - 6 \ge 0$,
this has $\det A = L(L-4)$  solutions  by {\sc Proposition} \ref{lem:basic}.
In this case there is no permutation redundancy to remove.
(iii) The equation (\ref{eq:sce3}) reads
\begin{equation*}
A^{1/2/3} \vec{u} =
\begin{pmatrix}
L-2 & 1 & 1 \\
1 & L-2 & 1 \\
1 & 1 & L-2
\end{pmatrix}
\begin{pmatrix} u_{1,1} \\ u_{1,2} \\ u_{1,3} \end{pmatrix}
\equiv
\vec{c}\quad \mathrm{mod}\ \Z^3
\end{equation*}
for some $\vec{c}$.
Here we have written $A$ as $A^{1/2/3}$ to match the notation
in Section \ref{subsec:general}.
Due to the assumption $P_1 = L - 6 \ge 0$, this has
$\det A^{1/2/3} = L(L-3)^2$ solutions in $(\R/\Z)^3$ by
{\sc Proposition} \ref{lem:basic}.
However they contain various diagonal solutions.
For example under the condition $u_{1,1} = u_{1,2}$,
the above equation  reduces to
\begin{equation*}
A^{12/3} \vec{u} =
\begin{pmatrix}
L-1 & 1 \\
2 & L-2
\end{pmatrix}
\begin{pmatrix} u_{1,1} \\ u_{1,3} \end{pmatrix}
\equiv
\vec{c}\,'\quad \mathrm{mod}\ \Z^2
\end{equation*}
for some $\vec{c }\,'$, which
has $\det A^{12/3} = L(L-3)$ solutions.
Similarly there are diagonal solutions counted by
$\det A^{13/2}$ and  $\det A^{23/1}$, which are both equal to $L(L-3)$.
Finally there is the completely diagonal one $u_{1,1} = u_{1,2} = u_{1,3}$
satisfying $A^{123} u_{1,1} = L u_{1,1} \equiv c''$ mod $\Z$
for some $c''$.
By the inclusion-exclusion principle, we can now compute the number of
off-diagonal solutions in (iii) as
\begin{equation}\label{eq:1112}
\begin{split}
&\frac{1}{3!}
\left(\det A^{1/2/3} - \det A^{12/3} - \det A^{13/2}
- \det A^{23/1} + 2\det A^{123}\right)\\
&\qquad = \frac{1}{6}L(L-4)(L-5),
\end{split}
\end{equation}
where we have removed the $3!$--fold $\frak{S}_3$ redundancy.
Assembling the contributions from (i), (ii) and (iii) we get
\begin{equation*}
L + L(L-4) + \frac{1}{6}L(L-4)(L-5) = \binom{L}{3}
\end{equation*}
as desired.

\subsection{General case}\label{subsec:general}
Let us proceed to the general case where both
the quantum space data $\nu = (\nu_s)$ and the magnon number
$M \in \Z_{\ge 0}$ are arbitrary.
The examples in Section \ref{subsec:examples} already elucidate
the essential feature in our counting.
We dare do some overcounting by 
classifying diagonal solutions in terms of their
patterns like $12/3$, and finally subtract them
via a kind of the inclusion-exclusion principle.
The most natural framework to systematize such a process is
the partition of sets and the M\"obius inversion trick
summarized in Appendix \ref{app:mobius}.
Compare the coefficients $1$, $-1$, $-1$,
$-1$, $2$ appearing in
(\ref{eq:1112}) with (\ref{eq:muexample}).
Below we will use the terminology and the notation therein.

We start with SCE (\ref{eq:sce3}).
Given $M \in \N$, fix a string pattern  $(N_m)$  satisfying
 (\ref{eq:mpartition}).
Set
\begin{equation}\label{eq:jdef}
{\mathcal J} = \{ j \in \N \mid N_j > 0\}, \quad j_0 = \max {\mathcal J}.
\end{equation}
Thus $d = \sum_{j \in \J}N_j$.
Take an element  $\pi = (\pi^{(1)}, \pi^{(2)}, \ldots, \pi^{(j_0)})
\in L_{N_1} \times \cdots \times L_{N_{j_0}}$ in the sense of
Appendix \ref{app:mobius}.
Thus  for each $k$,
$\pi^{(k)} = (\pi^{(k)}_1, \ldots,\pi^{(k)}_l)$ is a partition of the set
$\{1, \ldots, N_{k}\}$ into the blocks:
\begin{equation*}
\{1, \ldots, N_{k}\} = \pi^{(k)}_1 \sqcup \cdots  \sqcup \pi^{(k)}_l
\end{equation*}
for some $l$.
Here and in what follows, if $N_i = 0$ for $1 \le i < j_0$,
the corresponding $i$th component should be just dropped.
For example $L_{N_1} \times \cdots \times L_{N_{j_0}}$ actually means
the product of the set $L_{N_j}$ over $j \in \J$, and
$(\pi^{(1)}, \ldots,\pi^{({j_0})})$ is in fact
$(\pi^{(j)})_{j \in \J}$.
To classify the diagonal solutions, we introduce
\begin{displaymath}
\begin{split}
\Sol'_\pi = \{\vec{u} = (u_{k\beta})
 \mid u_{k\alpha} = u_{k\beta}
&\text{ if $\alpha$ and $\beta$ belong}\\
& \text{to the same block of $\pi^{(k)}$ for each $k$}\},\\
\Sol_\pi = \{\vec{u} = (u_{k\beta})
 \mid u_{k\alpha} = u_{k\beta}
& \text{ if and only if $\alpha$ and $\beta$ belong}\\
&\text{to the same
block of $\pi^{(k)}$ for each $k$}\}.
\end{split}
\end{displaymath}
By the definition it follows that ($\vert\, \cdot \, \vert$ here denotes
the cardinality.)
\begin{equation*}
\vert \Sol'_\pi \vert = \sum_{\pi' \le \pi} \vert \Sol_{\pi'} \vert
\end{equation*}
in terms of the partial order $\le$ of the poset
$L_{N_1} \times \cdots
 \times L_{N_{j_0}}$ introduced in Appendix \ref{subapp:5}.
By means of the M\"obius inversion formula
this is equivalent to
\begin{equation*}
\vert \Sol_\pi \vert = \sum_{\pi' \le \pi} \mu(\pi',\pi)
\vert \Sol'_{\pi'} \vert\qquad \text{for any } \pi \in
L_{N_1} \times \cdots \times L_{N_{j_0}}.
\end{equation*}
The off-diagonal solutions in question are counted  by setting
$\pi = \pi_{\text{max}}$, which is the maximal element of
$L_{N_1} \times \cdots \times L_{N_{j_0}}$
explained in Appendix \ref{subapp:5}.
Removing the $\frak{S}_{N_m}$ redundancy for $m \in \J$,
they are enumerated as
\begin{equation}
\frac{\vert \Sol_{\pi_{\text{max}}}\vert}{\prod_{m \in \J}N_m!}
= \frac{\sum_{\pi \in L_{N_1} \times \cdots \times L_{N_{j_0}}}
\mu(\pi,\pi_{\text{max}})
\vert \Sol'_{\pi}\vert}{\prod_{m \in \J}N_m!}.\label{eq:number}
\end{equation}
Let us  evaluate $\vert \Sol'_\pi \vert$.
In the original SCE (\ref{eq:sce3}),
impose the constraint
$u_{k\alpha} = u_{k\beta}$ on $\vec{u} = (u_{k\beta})$, if
$\alpha$ and $\beta$  belong to the same block of
$\pi^{(k)}$ for each $k \in \J$.
As we got $A^{12/3}$ from $A = A^{1/2/3}$ in Section
\ref{subsec:examples}, the result has the form:
\begin{equation}
A^\pi \vec{u}_\pi \equiv \vec{c}_\pi \quad 
\mathrm{mod}\
\Z^{l(\pi)}\qquad
\text{for  some}\ \vec{c}_\pi \in {\R}^{l(\pi)}.
\label{eq:api}
\end{equation}
Here $l(\pi)$ denotes the length of
$\pi = (\pi^{(1)}, \ldots, \pi^{({j_0})})$ as specified in
Appendix \ref{subapp:5}.
In the new unknown  $\vec{u}_\pi = (u_{k\beta})$,
$\beta$ is now labeled by the blocks of $\pi^{(k)}$ for each
$k \in  \J$.
$A^\pi$ is an $l(\pi)$ by  $l(\pi)$ integer matrix obtained by
a reduction of the matrix $A$  which was $d$ by $d$ originally
($d = \sum_{j \in \J}N_j$).
It is formed by summing up the $(k\beta)$ columns of $A$
over those $\beta$ belonging to the same block of $\pi^{(k)}$,
and discarding all but one rows for each  block.
For example, suppose $\J = \{1,2\}$ .
Thus only $N_1$ and $N_2$ are non-zero and
$\pi = (\pi^{(1)}, \pi^{(2)}) \in L_{N_1} \times L_{N_2}$.
Let
$\pi^{(1)} = (\pi^{(1)}_1, \pi^{(1)}_2, \pi^{(1)}_3)$ and
$\pi^{(2)} = (\pi^{(2)}_1, \pi^{(2)}_2)$ so that
$l(\pi^{(1)}) = 3, l(\pi^{(2)}) = 2$ and $l(\pi) = 5$.
Denote the number of elements in the block  $\pi^{(1)}_i$ by $\lambda_i$
(resp.\ $\vert \pi^{(2)}_i\vert$ by $\mu_i$).
By the definition $\lambda_1 + \lambda_2 + \lambda_3 = N_1, \,
\mu_1 + \mu_2 = N_2$, and
the matrix $A^\pi$ reads
(${\mathcal P}_i = P_i + N_i$)
\begin{equation*}
\left(
\begin{array}{ccccc}
{\mathcal P}_1 + \lambda_1 & \lambda_2 & \lambda_3 & 2\mu_1 & 2\mu_2 \\
\lambda_1 & {\mathcal P}_1 + \lambda_2 & \lambda_3 & 2\mu_1 & 2\mu_2 \\
\lambda_1 & \lambda_2 & {\mathcal P}_1 + \lambda_3  & 2\mu_1 & 2\mu_2 \\
2\lambda_1 & 2\lambda_2 & 2\lambda_3 & {\mathcal P}_2 + 3\mu_1 & 3\mu_2 \\
2\lambda_1 & 2\lambda_2 & 2\lambda_3 & 3\mu_1 & {\mathcal P}_2 + 3\mu_2
\end{array}
\right).
\end{equation*}
This is easily seen from the example of $A$ in Section \ref{subsec:rule}.
In general $A^\pi$ is not symmetric.
Its matrix element is given by
$A^\pi_{(m,i), (k,j)} = \delta_{m,k}\delta_{i,j}(P_m+N_m) +
(2\min(m,k)-\delta_{m,k})\vert \pi^{(k)}_j \vert$ for
$m, k \in \J$,
$1 \le i \le l(\pi^{(m)})$ and
$1 \le j \le l(\pi^{(k)})$.
Note that $A = A^{\pi_{\text{max}}}$.

By the definition $\vert \Sol'_\pi \vert$ counts the number of all the
solutions to (\ref{eq:api}).
Therefore from {\sc Proposition} \ref{lem:basic} we have
\begin{equation}
\vert \Sol'_\pi \vert = \vert \det A^\pi \vert.
\label{eq:solapi}
\end{equation}
It is straightforward to show
\begin{lemma}\label{lem:detapi}
For any $\pi \in L_{N_1} \times \cdots \times L_{N_{j_0}}$ we have
\begin{align}
\det A^\pi &= \det_{m,k \in \J}(F_{m,k})
\prod_{m \in \J} (P_m + N_m)^{l(\pi^{(m)})-1},\label{eq:detapi}\\
F_{m,k} &= \delta_{m,k}P_m + 2\min(m,k)N_k.
\label{eq:fdef}
\end{align}
\end{lemma}
The lemma holds without assuming (\ref{eq:ppositive}).
The dependence on $\pi$ enters only through $l(\pi^{(m)})$.
Like $A^{\pi}$, the matrix $F$ is non-symmetric in general:
\begin{displaymath}
F =
\left(
\begin{array}{ccccc}
P_1+2N_1 & 2N_2 & 2N_3 & \cdots & 2N_{j_0} \\
2N_1 & P_2+4N_2 & 4N_3 & \cdots & 4N_{j_0} \\
2N_1 & 4N_2 & P_3+6N_3 & \cdots & \vdots \\
\vdots & \vdots & \vdots & \ddots & 2({j_0}-1)N_{j_0} \\
2N_1 & 4N_2 & 6N_3 & \cdots & P_{j_0}+2{j_0}N_{j_0}
\end{array}
\right).
\end{displaymath}
We remark that 
\begin{equation}\label{eq:fsum}
\sum_{k \ge 1} F_{m,k} = \gamma_m
\end{equation}
for any $m$, where $\gamma_m$ is defined in (\ref{eq:gammadef}).
\begin{lemma}\label{lem:positive}
If  $P_m \ge 0$ for any $m \in \J$, then $\det A^\pi > 0$.
\end{lemma}
\begin{proof}
Denote $\det_{m,k \in \J}(F_{m,k})$ simply by $\det_{\J} F$.
By {\sc Lemma} \ref{lem:detapi}, it suffices to verify
$\det_{\J}F > 0$.
We do this by a double induction
on $\vert {\mathcal J}\vert$ and $\sum_m P_m$
regarding $P_m$ as non-negative variables independent of $\{N_m \}$.
First let ${\mathcal J} = \{j_1 < \cdots
 < j_l\}$ be arbitrary and $\sum_m P_m=0$.
Thus $P_m = 0$ for any $m \in {\mathcal J}$, hence
$\det_{\J}F = (\prod_{j \in {\mathcal J}}2N_j) \det_{m,k \in {\mathcal
J}}(\min(m,k))
= (\prod_{j \in {\mathcal J}}2N_j)j_1(j_2-j_1)\cdots
(j_l-j_{l-1})>0$.
Next let ${\mathcal J} = \{j\}$.
Then $\det_{\J}F = P_j + 2N_j > 0$ because of the assumption $P_j \ge 0$.
Finally let ${\mathcal J}$ and $\sum_mP_m > 0$ be arbitrary.
Then there exists $i \in {\mathcal J}$ such that $P_i > 0$.
Setting $P'_j = P_j - \delta_{j,i}$ and
${\mathcal J}' = {\mathcal J} \setminus  \{i\}$,
one can expand the determinant as
$\det_{\mathcal J}F(\{P_m\}) =
\det_{{\mathcal J}'}F(\{P_m\}) + \det_{\mathcal J}F(\{P'_m\})$.
By induction the two terms in the right hand side are both positive.
\end{proof}
\par
Although it is more direct to expand $\det_{m,k \in \J}(F_{m,k})$ from the
beginning,
we have presented a proof in the above form because it generalizes
to arbitrary simple Lie algebra case that will be treated in our
subsequent paper.

The specialization $\pi = \pi_{\text{max}}$ in the above leads to
\begin{corollary}\label{cor:detanonzero}
If $P_m \ge 0$ for any $m$, then $\det A > 0$.
\end{corollary}
Given the quantum space data $\nu = (\nu_s)$ and the
string pattern $N = (N_m)$, we define
\begin{equation}
R(\nu,N) =
\det_{m,k \in {\mathcal J}}(F_{m,k})
\prod_{m \in \J} \frac{1}{N_m}
\binom{P_m + N_m - 1}{N_m - 1},
\label{eq:Rdef}
\end{equation}
when $\J \neq \emptyset$.
Here $P_m = P_m(\nu,N),\; F_{m,k}$ and $\J$ are given by
(\ref{eq:pdef}), (\ref{eq:fdef})  and (\ref{eq:jdef}), respectively.
When $\J = \emptyset$, namely $\forall N_m = 0$, we set
$R(\nu,N) = 1$ irrespective of $\nu$.
In the definition itself, we do not need to assume that
$P_m \ge 0$ for those $m \in \J$, and $(\nu_s)$ can be arbitrary
complex parameters.
\begin{theorem}\label{th:Rnds}
Assume that $P_m \ge 0$ for any  $m \in \J$.
Then the number of off-diagonal solutions to
SCE (\ref{eq:sce3}) is equal to $R(\nu,N)$.
\end{theorem}
\begin{proof}
Under the assumption the number of off-diagonal solutions has
already been obtained in (\ref{eq:number}).
By virtue of (\ref{eq:solapi}) and {\sc Lemma} \ref{lem:positive}
it is equal to
\begin{displaymath}
\begin{split}
&\frac{1}{\prod_{m \in \J}N_m!}
\sum_{\pi \in L_{N_1} \times \cdots \times L_{N_{j_0}}}
\mu(\pi,\pi_{\text{max}})
\det A^\pi \\
&\qquad= \det_{m,k \in {\mathcal J}}(F_{m,k})
\sum_{\pi \in L_{N_1} \times \cdots
 \times L_{N_{j_0}}} \mu(\pi,\pi_{\text{max}})
\prod_{m \in {\mathcal J}}
\frac{(P_m+N_m)^{l(\pi^{(m)})-1}}{N_m!},\\
\end{split}
\end{displaymath}
where we have substituted {\sc Lemma} \ref{lem:detapi}.
By means of (\ref{eq:musum})  the $\pi$--sum can be taken, leading to
\begin{displaymath}
\det_{m,k \in {\mathcal J}}(F_{m,k})
\prod_{m \in {\mathcal J}}
\frac{(P_m+N_m)_{N_m}}{N_m!(P_m+N_m)} = R(\nu,N).
\end{displaymath}
\end{proof}
\par
Under the specialization $\nu_s = L\delta_{1,s}$,
the above $R(\nu,N)$ reproduces the number of off-diagonal
solutions for each string pattern exemplified in Section \ref{subsec:examples}.

Expanding the determinant in (\ref{eq:Rdef}), one can rewrite $R(\nu,N)$
as follows:
\begin{align}
R(\nu,N) &= \sum_{J \subset \N} D_J
\prod_{m \in \N \setminus J}\binom{P_m + N_m}{N_m}
\prod_{m \in J}\binom{P_m + N_m - 1}{N_m - 1},\label{eq:Rexpand}\\
D_J &= \begin{cases}
1 & \text{ if } J = \emptyset,\\
\det_{m,k \in J}(2\min(m,k)-\delta_{m,k}) & \text{otherwise}.
\end{cases}\label{eq:Ddef}
\end{align}
In deriving this we have used
$\binom{\;\ast \;}{0} = 1,\; \binom{\;\ast\;}{-1} = 0$.
{}From this expression $R(\nu,N) \in \Z$ is manifest if
$\forall \nu_s \in \Z$.
In Section \ref{sec:multiplicity} we will mainly work with the formula
(\ref{eq:Rexpand}) rather than (\ref{eq:Rdef}).


\section{$R(\nu,N)$ as weight multiplicity}\label{sec:multiplicity}

Set
\begin{equation*}
K(\nu,N) = \prod_{m \ge 1}\binom{P_m + N_m}{N_m}.
\end{equation*}
This is a generalization of Bethe's fermionic formula in 
(\ref{eq:completeness0}) corresponding to 
the quantum space (\ref{eq:qspace})  \cite{K}.
The nature of our $R(\nu,N)$ becomes most transparent
by a parallel analysis on $K(\nu,N)$.
It contains $K(\nu,N)$ as the summand in (\ref{eq:Rexpand})
corresponding to $J = \emptyset$.

In this section we fix $l \in \N$.
It plays a role of ``cut-off\/" similar to  $j_0$ in
(\ref{eq:jdef}) and has nothing to do with the length function of
partitions.
We will introduce various functions indexed with $l$, which
tend to the quantities in our problem in the limit
$l \rightarrow \infty$.
In particular,  $P_m$  and $\gamma_m$
in Sections \ref{subsec:RK}--\ref{subsec:analytic}
stand for the truncations of
(\ref{eq:pdef})--(\ref{eq:gammadef}) by $l$:
\begin{align}
P_m &= P_m(\nu,N) = \gamma_m - 2\sum_{k = 1}^l \min(m,k)N_k,
\label{eq:truncatedp}\\
\gamma_m &= \gamma_m(\nu) = \sum_{k =1}^l \min(m,k) \nu_k.
\label{eq:truncatedgamma}
\end{align}
We do not prepare new symbols for them as they will be used only
in the said sections.
We set $\N_l = \{1,2,\ldots, l\}$.
The binomial coefficient is the one specified in the beginning of
Section \ref{subsec:examples}.

\subsection{$R_l(\nu,N)$ and $K_l(\nu,N)$}\label{subsec:RK}
Let $\nu = (\nu_s)$, $\nu_1, \ldots,  \nu_l  \in \C$ and
$N = (N_m)$,  $N_1$, \dots, $N_l \in \Z_{\ge 0}$ be arbitrary.
Define
\begin{align}
R_l(\nu,N) &= \sum_{J \subset \N_l} D_J
\prod_{m \in \N_l\setminus J}\binom{P_m + N_m}{N_m}
\prod_{m \in J}\binom{P_m + N_m - 1}{N_m - 1},\label{eq:Rldef}\\
K_l(\nu,N) &=
\prod_{m \in \N_l}\binom{P_m + N_m}{N_m},\label{eq:Kldef}
\end{align}
where $D_J$ is specified by (\ref{eq:Ddef}).
When $N = 0$ (i.e., $\forall N_m = 0$), we have
$R_l(\nu,0) = K_l(\nu,0) = 1$ irrespective of $\nu$.
Obviously one has
$R(\nu,N) = \lim_{l \rightarrow \infty}R_l(\nu,N)$ and
$K(\nu,N) = \lim_{l \rightarrow \infty}K_l(\nu,N)$, where the
limits render no subtlety.
We will utilize two other expressions of $R_l(\nu,N)$.
The first one is the analogue of (\ref{eq:Rdef}):
\begin{equation}
R_l(\nu,N) =
\left(\det_{m,k \in \{i \in \N_l \mid N_i \neq 0\}}F_{m,k}\right)
\prod_{m \in \N_l, N_m \neq 0} \frac{1}{N_m}
\binom{P_m + N_m - 1}{N_m - 1}.
\label{eq:Rldef2}
\end{equation}
To match (\ref{eq:Rldef}), the right side of this
should be understood as 1 when $N = 0$.
To deduce the second expression, for $J \subset \N_l$  we introduce
\begin{alignat}{2}
\nu[J] &= (\nu[J]_s),
& \quad \nu[J]_s &= \nu_s - 2\theta(s \in J),
\label{eq:nuJ}\\
N[J]  &= (N[J]_m),
& \quad N[J]_m &= N_m - \theta(m \in J),
\label{eq:nJ}
\end{alignat}
where $\theta(\text{true}) = 1$ and
$\theta(\text{false}) = 0$.
With the aid of
\begin{equation*}
P_m(\nu[J],N[J]) = P_m(\nu,N),
\end{equation*}
one can rewrite (\ref{eq:Rldef}) as
\begin{equation}
R_l(\nu,N) = \sum_{J \subset \N_l} D_J
\prod_{m \in \N_l}\binom{P[J]_m + N[J]_m}{N[J]_m}.\label{eq:Rldef3}
\end{equation}

\subsection{Generating functions}\label{subsec:generating}

Let us introduce generating functions:
\begin{align}
R_l(\nu \vert w) &= \sum_N R_l(\nu,N)w_1^{N_1} \cdots
w_l^{N_l},\label{eq:rlwdef}\\
K_l(\nu \vert w) &= \sum_N K_l(\nu,N)w_1^{N_1} \cdots
w_l^{N_l},\label{eq:klwdef}
\end{align}
where $w = (w_1, \ldots, w_l)$ and
$\sum_N$ extends over $N_1$, \dots, $N_l \in \Z_{\ge 0}$.
\begin{proposition}\label{pr:nuzero}
When $\nu = 0$ (i.e., $\forall \nu_s = 0$) we have
\begin{equation*}
R_l(0 \vert w) = 1.
\end{equation*}
\end{proposition}
\begin{proof}
We show $R_l(0, N) = 0$ for any $N \neq 0$.
Note that $\nu = 0$ implies $\forall \gamma_m = 0$.
Therefore the assertion follows from the expression (\ref{eq:Rldef2})
and (\ref{eq:fsum}).
\end{proof}
This simple observation  will eventually lead to a non-trivial
consequence (\ref{eq:FOP}) whose derivation 
is analogous to the ``denominator formula".
In contrast, $K_l(0 \vert w)$ is not a simple function.
See (\ref{eq:Kexp}).
\subsection{Analytic formula for generating functions}\label{subsec:analytic}

Consider the variables
$\{z_{j,i-1} \mid  1 \le i \le  j \le l \}$
related via
\begin{equation}\label{eq:zsystem}
z_{j,i} = z_{j,i-1}
\left(1-z_{i, i-1} \right)^{-2(j-i)}\qquad
1\leq i < j \leq l.
\end{equation}
For $1 \le i \le l$ we define the function
$\psi_i$ by
\begin{equation}\label{eq:psidef}
\psi_i = \prod_{j=i}^l
\left(1-z_{j, j-1}\right)^{-\beta_j - 1} \qquad 1 \le i \le l,
\end{equation}
where $\beta_1$, \dots, $\beta_l \in \C$ are  parameters.
\begin{lemma}\label{lem:hkoty}
$\psi_i$ has a formal power series  expansion
\begin{displaymath}
\psi_i = \sum_{\{N_j\}} \prod_{j=i}^l
\binom{\beta_j + 2\sum_{j < k \le l}(k-j)N_k + N_j}{N_j}
\left(z_{j,i-1}\right)^{N_j} \qquad 1 \le i \le l,
\end{displaymath}
where the sum $\sum_{\{N_j\}}$ extends over
$N_i$, \dots, $N_l \in \Z_{\ge 0}$.
\end{lemma}
\begin{proof}
We prove by induction on $i$.
The case $i=l$ is due to the formula
\begin{equation}\label{eq:exp}
(1-z)^{-\beta-1} = \sum_{N=0}^\infty \binom{\beta+N}{N}z^N.
\end{equation}
Assume $\psi_{i+1}$ has the above expansion.
Then from (\ref{eq:psidef}) $\psi_i$ is 
\begin{displaymath}
 (1-z_{i,i-1})^{-\beta_i - 1}
\sideset{}{'}\sum_{ \{N_j\} }
\prod_{j=i+1}^l
\binom{\beta_j + 2\sum_{j < k \le l}(k-j)N_k + N_j}{N_j}
\left(z_{j,i}\right)^{N_j},
\end{displaymath}
where the sum $\sum'_{\{N_j\}}$ is over
$N_{i+1}$, \dots, $N_l \in \Z_{\ge 0}$.
Upon substituting (\ref{eq:zsystem}), the right hand side becomes
\begin{displaymath}
\begin{split}
\sideset{}{'}\sum_{\{N_j\}}
\biggr\{ &
(1-z_{i,i-1})^{-\beta_i -2\sum_{i < k \le l}(k-i)N_k -1}\\
&\qquad
\times \prod_{j=i+1}^l\binom{\beta_j +
2\sum_{j < k \le l}(k-j)N_k + N_j}{N_j}
\left(z_{j,i-1}\right)^{N_j}\biggl\}.
\end{split}
\end{displaymath}
Applying (\ref{eq:exp}) again,
we obtain the desired expansion.
\end{proof}
\par
This lemma is originally due to \cite{K}.
Here we have quoted the version reproduced in \cite{HKOTY}.

In the sequel,
we will only work with the variables
$z_j = z_{j,0},\, v_j = z_{j,j-1}$ and $w_j$
($1 \le j \le l$):
\begin{equation}\label{eq:zvwrel}
v_j = z_j \prod_{k=1}^{j-1}(1-v_k)^{-2(j-k)}, \quad
w_j = z_j \prod_{k=1}^l(1-v_k)^{-2j} \quad
\end{equation}
where the former relation is due to  (\ref{eq:zsystem}), while the
latter is  the definition of $w_j$.
($z_i$ here should not be confused with 
the $z_i(q)$ in Section \ref{sec:beq0}.)
Note that
$v_i = w_i\prod_{k=1}^l(1-v_k)^{2\min(i,k)}$.
Let $z = (z_1,\ldots, z_l)$ and
$w= (w_1,\ldots, w_l)$.
We will denote $z_1 = \cdots = z_l = 0$ simply by
$z = 0$, and
${dz_1 \wedge \ldots \wedge dz_l}/{z_1 \cdots
 z_l}$ by
${dz}/{z}$, etc.
The variables $z$ are holomorphic functions of
$w$ around $w = 0$.
This is due to $w = 0$ and
${\partial w_i}/{\partial z_j} = \delta_{i j}$
at $z=0$.
Setting $i=1$ in {\sc Lemma} \ref{lem:hkoty}, we have
\begin{equation}
\label{eq:basic}
\prod_{j=1}^l
\binom{\beta_j + 2\sum_{j < k \le l}(k-j)N_k + N_j}{N_j}
= \Res_{z = 0}\left(
\prod_{j=1}^l(1-v_j)^{-\beta_j-1}z_j^{-N_j} \right)
\frac{dz}{z}.
\end{equation}
Under a further specialization to
$\beta_j = \gamma_j - 2\sum_{k=1}^lkN_k$, this becomes
\begin{align}
\prod_{m \in \N_l}\binom{P_m + N_m}{N_m} &=
\Res_{z=0}
\left(
\prod_{j=1}^l(1-v_j)^{-\gamma_j-1}w_j^{-N_j}\right)
\frac{dz}{z}\notag \\
&= \Res_{w=0}
\left(
\prod_{j=1}^l(1-v_j)^{-\gamma_j-l(l+1)-1}w_j^{-N_j}\right)
\frac{\partial z}{\partial w}
\frac{dw}{w},\label{eq:resK}
\end{align}
where ${\partial z}/{\partial w}$
represents the Jacobian $\det_{i,j \in \N_l}
({\partial z_i}/{\partial w_j})$.
In (\ref{eq:resK}) replace $N_m$ by $N[J]_m$, $P_m$ by $P_m(\nu[J],N[J])$
defined in (\ref{eq:nuJ})--(\ref{eq:nJ}),
and $\gamma_j$ by $\gamma_j - 2\sum_{i \in J}\min(j,i)$.
The result reads
\begin{equation}
\begin{split}
&\prod_{m \in \N_l}\binom{P[J]_m + N[J]_m}{N[J]_m}\\
&\qquad
 =
\Res_{w=0}
\left(
\prod_{j=1}^l
(1-v_j)^{-\gamma_j-l(l+1)-1}w_j^{-N_j}\right)
\left(\prod_{i \in J}v_i\right)
\frac{\partial z}{\partial w}
\frac{dw}{w}.\label{eq:resR}
\end{split}
\end{equation}
Notice that the left side constitutes the summand in (\ref{eq:Rldef3}).
{}In terms of the generating functions the results
(\ref{eq:resK})--(\ref{eq:resR}) are
stated as
\begin{proposition}\label{pr:RKgenerating}
\begin{align}
K_l(\nu\vert w) &= K_l(0\vert w)
 \prod_{j=1}^l(1-v_j)^{-\gamma_j},\\
K_l(0\vert w) &= \frac{\partial z}{\partial w}
\prod_{j=1}^l(1-v_j)^{-l(l+1)-1},
\label{eq:Kexp}
\\
R_l(\nu\vert w) &= R_l(0\vert w)
\prod_{j=1}^l(1-v_j)^{-\gamma_j},
\\
R_l(0\vert w) &= K_l(0\vert w)\left(\sum_{J \subset \N_l}D_J
\prod_{i \in J}v_i\right).
\end{align}
\end{proposition}
Combining this with {\sc Proposition} \ref{pr:nuzero},
we find
\begin{theorem}\label{th:mainl}
\begin{gather}
R_l(\nu\vert w) = \frac{K_l(\nu\vert w)}{K_l(0\vert w)}
= \prod_{j=1}^l(1-v_j)^{-\gamma_j},\label{eq:RKK}\\
R_l(\nu\vert w)R_l(\nu'\vert w) = R_l(\nu+\nu'\vert
w),\label{eq:factorization}\\
\frac{\partial w}{\partial z}
\prod_{j=1}^l(1-v_j)^{l(l+1)+1}
= \sum_{J \subset \N_l}D_J\prod_{i \in J}v_i,\label{eq:FOP}
\end{gather}
where $\nu + \nu' = (\nu_s + \nu'_s)_{s=1}^l$ and
$\gamma_j$ is the truncated one (\ref{eq:truncatedgamma}).
\end{theorem}
The factorization property  (\ref{eq:factorization}) is
enjoyed only by $R_l$  and  not by $K_l$.
It is due to (\ref{eq:RKK}) and
$\gamma(\nu+\nu')_j = \gamma(\nu)_j + \gamma(\nu')_j$.

\subsection{$l \rightarrow \infty$ limit}\label{subsec:linfty}

Let
$R(\nu\vert w) = \lim_{l \rightarrow \infty}R_l(\nu\vert w)$
and
$K(\nu\vert w) = \lim_{l \rightarrow \infty}K_l(\nu\vert w)$ be
formal power series in infinitely many variables
$w=(w_j)_{j\ge 1}$.
They can also be viewed as the series in  $(v_j)_{j\ge 1}$ upon the
substitution
$w_i = v_i\prod_{k\ge 1}(1-v_k)^{-2\min(i,k)}$.
See the remark after (\ref{eq:zvwrel}).
In the $l \rightarrow \infty$ limit {\sc Theorem} \ref{th:mainl}
yields
\begin{theorem}\label{th:main}
\begin{gather}
R(\nu\vert w) = \frac{K(\nu\vert w)}{K(0\vert w)}
= \prod_{j \ge 1}(1-v_j)^{-\gamma_j},\label{eq:RKKinf}\\
R(\nu\vert w)R(\nu'\vert w) = R(\nu+\nu'\vert w),\label{eq:factorizationinf}
\end{gather}
where $\nu + \nu' = (\nu_s+\nu'_s)_{s \ge 1}$ and
$\gamma_j$ is defined by (\ref{eq:gammadef}).
\end{theorem}
We specialize $R(\nu\vert w)$ and $K(\nu\vert w)$ as follows:
\begin{align}
R(\nu) &:= e^{\gamma_\infty(\nu)\Lambda_1}R(\nu\vert w)\vert_{w_j =
e^{-j\alpha_1}}
= \sum_N R(\nu,N)x^{\gamma_\infty(\nu)-2\sum_{j\ge 1}jN_j},
\label{eq:Rspecial}\\
K(\nu) &:= e^{\gamma_\infty(\nu)\Lambda_1}K(\nu\vert w)\vert_{w_j =
e^{-j\alpha_1}}
= \sum_N K(\nu,N)x^{\gamma_\infty(\nu)-2\sum_{j\ge 1}jN_j},
\label{eq:Kspecial}
\end{align}
where the sum $\sum_N$ runs over
$N_1$, $N_2$, $\ldots \in \Z_{\ge 0}$,
$\alpha_1$ and $\Lambda_1$ are the simple root and the fundamental weight,
respectively.
$x = e^{\Lambda_1}$ is a formal variable and
$\gamma_\infty(\nu) = \sum_{s\ge 1}s\nu_s$ in accordance with
(\ref{eq:gammadef}).
{}From {\sc Theorem} \ref{th:main} it follows that
\begin{align}
R(\nu) &= \frac{K(\nu)}{K(0)},\label{eq:RKK2}\\
R(\nu)R(\nu') &= R(\nu+\nu').\label{eq:factorizationspecial}
\end{align}
We remark that {\sc Proposition} \ref{pr:RKgenerating},
{\sc Theorem} \ref{th:mainl} and \ref{th:main} are all valid
for $\nu_s \in \C$.
The specialization $w_j = e^{-j\alpha_1}$
induces an effect also for $v_j$ (\ref{eq:zvwrel}) and
$\prod_{j \ge 1}(1-v_j)^{-\gamma_j}$ in (\ref{eq:RKKinf}).
It will be worked out in Section \ref{subsec:qsys}.

\subsection{Combinatorial completeness}
\label{subsec:qsys}
{}From now on we assume that $\forall \nu_s \in \Z_{\ge 0}$.
For $m \in \Z_{\ge 0}$ we define
\begin{align}
\delta_m &= (\nu_s), \quad \nu_s = \delta_{s, m},\notag\\
Q_m &= R(\delta_m) \in x^m\C[[ x^{-2}]].\label{eq:Qdef}
\end{align}
{}From the decomposition
$\nu = (\nu_s) = \sum_{s\ge 1}\nu_s\delta_s$
and (\ref{eq:factorizationspecial}) we have
\begin{equation}\label{eq:completeness1}
R(\nu) = \prod_{s \ge 1} Q_s^{\nu_s}
\end{equation}
for general $\nu = (\nu_s)$.
\begin{proposition}\label{pr:sum}
\begin{equation*}
R(\lambda) = R(\mu) + R(\nu),
\end{equation*}
where $\lambda = (\lambda_s), \mu = (\mu_s)$ and $\nu = (\nu_s)$ are
related as ($s \in \N$)
\begin{displaymath}
\lambda_s = \nu_s + 2\delta_{s, k},\quad
\mu_s = \nu_s + \delta_{s, k+1} + \delta_{s, k-1}\qquad
\text{for some $k \in \N$}.
\end{displaymath}
\end{proposition}
\begin{proof}
Put $N' = (N'_m)$, $N'_m = N_m - \delta_{m,k}$.
Then it is easy to check
\begin{gather}
\gamma_\infty(\lambda) - 2\sum_{j \ge 1} jN_j =
\gamma_\infty(\mu) - 2\sum_{j \ge 1} jN_j =
\gamma_\infty(\nu) - 2\sum_{j \ge 1} jN'_j,\label{eq:powers}\\
P_m(\lambda,N) = P_m(\nu,N') = P_m(\mu,N) + \delta_{m,k}.\label{eq:ps}
\end{gather}
{}From (\ref{eq:Rspecial}) and (\ref{eq:powers}) we are to show
$R(\lambda,N) = R(\mu,N) + R(\nu,N')$.
By expanding a binomial coefficient in
(\ref{eq:Rexpand}), the $R(\lambda,N)$ is expressed as
($P_m = P_m(\lambda,N)$)
\begin{displaymath}
\begin{split}
&\sum_{J \subset \N, \, k \not\in J} D_J
\left( A + B \right)
\prod_{m \in \N \setminus J,\, m \neq k}\binom{P_m + N_m}{N_m}
\prod_{m \in J}\binom{P_m + N_m - 1}{N_m - 1}\\
+ &\sum_{J \subset \N, \, k \in J} D_J
\prod_{m \in \N \setminus J}\binom{P_m + N_m}{N_m}
\left( C + D \right)
\prod_{m \in J, \, m \neq k}\binom{P_m + N_m - 1}{N_m - 1},
\end{split}
\end{displaymath}
where
$A = \binom{P_k + N_k - 1}{N_k}$,
$B = \binom{P_k + N_k - 1}{N_k-1}$,
$C = \binom{P_k + N_k - 2}{N_k-1}$ and
$D = \binom{P_k + N_k - 2}{N_k-2}$.
By  (\ref{eq:ps}) they can also be written as
$A = \binom{P_k(\mu,N) + N_k}{N_k}$,
$B = \binom{P_k(\nu,N')+N'_k}{N'_k}$,
$C = \binom{P_k(\mu,N) + N_k - 1}{N_k - 1}$ and
$D = \binom{P_k(\nu,N') + N'_k - 1}{N'_k - 1}$.
Thus the contributions containing $A$ and $C$ (resp.\
 $B$ and $D$)
amount to $R(\mu,N)$ (resp.\ $R(\nu,N')$).
\end{proof}
\par
Let $\mathbb{Q}((x))$ denote the field
of the formal Laurent series in $x$ over $\mathbb{Q}$
with
finitely many negative powers.
Clearly, $Q_m\in \mathbb{Q}((x^{-1}))$.
\begin{proposition}\label{pr:qsys}
%
(i) $Q_m$ satisfies
\begin{itemize}
\item[\textit{\/(a)}]
(recursion relation)
\begin{equation*}
Q_0 = 1, \quad Q_k^2 = Q_{k+1}Q_{k-1} + 1\qquad k \in \N,
\end{equation*}
\item[\textit{(b)}]
(asymptotic property)
\begin{equation*}
\lim_{k \rightarrow \infty}\frac{Q_{k+1}}{Q_k} = x.
\end{equation*}
\end{itemize}
\par
(ii) Conversely, the properties (a) and (b) above
 characterize
the series $Q_m \in \mathbb{Q}((x^{-1}))$.
\end{proposition}
\begin{proof}
(i). (a) Put $\nu = 0$ in {\sc Proposition} \ref{pr:sum} and apply
(\ref{eq:completeness1}).
(b) It is enough to show that the limit
$\lim_{k\to \infty} x^{-k}Q_k$ exists in $\mathbb{Q}
[[x^{-1}]]$.
Note that $P_m(\delta_k,N) = P_m(\delta_{k+1},N) - \theta(m \ge
k+1)$ {}from (\ref{eq:pdef}).
In the series $x^{-k}Q_k = x^{-k}R(\delta_k)$ in 
(\ref{eq:Rspecial}), those $N = (N_m)$ containing 
$N_j >0$ with  $j \ge k+1$  make contributions 
in the order higher than $2k+1$.
It follows that
$x^{-k}Q_k\equiv x^{-k-1}Q_{k+1}$ mod $x^{-2k-2}
\mathbb{Q}[[x^{-1}]]$.
Then, we have
\begin{displaymath}
x^{-k}Q_k\equiv x^{-k-1}Q_{k+1}\equiv
x^{-k-2}Q_{k+2}\equiv \cdots
\quad \text{mod}\ x^{-2k-2}\mathbb{Q}[[x^{-1}]],
\end{displaymath}
which means $\lim_{k\to \infty}x^{-k}Q_k$ exists.
%
%
%
%
%
%
%
%
%
(ii). Suppose $\tilde{Q}_m$ satisfies (a) and (b).
Setting $v_j = 1-\frac{\tilde{Q}_{j-1}\tilde{Q}_{j+1}}{\tilde{Q}^2_j}$,
we find
\begin{equation*}
\prod_{j=1}^l(1-v_j)^{-\gamma_j} =
\left(\frac{\tilde{Q}_l}{\tilde{Q}_{l+1}}\right)^{\gamma_l}
\prod_{j=1}^l\tilde{Q}^{\nu_j}_j,\quad
w_j = \left(\frac{\tilde{Q}_l}{\tilde{Q}_{l+1}}\right)^{2j}
\end{equation*}
by (\ref{eq:zvwrel}).
($\gamma_j$ here is the truncated one (\ref{eq:truncatedgamma}).)
Therefore (\ref{eq:RKK}) specializes to
\begin{equation*}
\sum_N R_l(\nu,N) \prod_{j=1}^l
\left(\frac{\tilde{Q}_l}{\tilde{Q}_{l+1}}\right)^{2jN_j} =
\left(\frac{\tilde{Q}_l}{\tilde{Q}_{l+1}}\right)^{\gamma_l}
\prod_{j=1}^l \tilde{Q}^{\nu_j}_j,
\end{equation*}
where $\sum_N$ is over $N_1$, \dots, $N_l \in \Z_{\ge 0}$.
By taking the limit $l \rightarrow \infty$ using (b)
for $\tilde{Q}_m$, this leads to
$R(\nu) = \prod_{j \ge 1} \tilde{Q}^{\nu_j}_j$.
Since $\nu_j$'s are arbitrary, we obtain
$\tilde{Q}_m = R(\delta_m)$.
Comparing this with (\ref{eq:Qdef})  we conclude
$\tilde{Q}_m = Q_m$.
\end{proof}
\par
It is immediate to check that the character of
the $(m+1)$-dimensional
irreducible $U_q(\selh)$-module $W_m$
(character with respect to the classical Cartan subalgebra)
\begin{equation*}
\mathrm{ch}\,W_m = \frac{x^{m+1}-x^{-m-1}}{x-x^{-1}},
\quad x =
e^{\Lambda_1}
\end{equation*}
fulfills the properties (a) and (b) in {\sc Proposition} \ref{pr:qsys}.
Thus from  (ii) we have
\begin{proposition}\label{cor:Qcharacter}
\begin{equation*}
Q_m = \mathrm{ch}\, W_m
 \qquad m \in \Z_{\ge 0}.
\end{equation*}
\end{proposition}
Our main result in Section \ref{sec:multiplicity} is 
the following.
\begin{theorem}[Combinatorial completeness]\label{th:completeness3}
Let 
$W(\nu)$
be the quantum space in (\ref{eq:qspace}),
$W(\nu) = \bigotimes_{s\geq 1} (W_s)^{\otimes \nu_s}$.
\par
(i) 
\begin{align*}
R(\nu) &= \mathrm{ch}\,W(\nu),\\
\sum_{N}{}^{(\lambda)} R(\nu,N) &= \dim W(\nu)_\lambda
\qquad \lambda \in \Z \Lambda_1.
\end{align*}
Here the sum $\sum_{N}^{(\lambda)}$ extends over
$N_1$, $N_2$, $\dots \in \Z_{\ge 0}$ such that
$\sum_{j \ge 1}j(\nu_j - 2N_j)\Lambda_1 = \lambda$, and 
$\dim W(\nu)_\lambda$ denotes the multiplicity of the weight $\lambda$.
In particular $R(\nu)$ is invariant under the Weyl group.
\par
(ii) (Kirillov \cite{K})
\begin{align*}
K(\nu) &= (1-e^{-\alpha_1})\, \mathrm{ch}\,W(\nu), \\
\sum_{N}{}^{(\lambda)}
 K(\nu,N) &= [W(\nu): V_\lambda] \ 
\qquad\lambda \in (\Z_{\ge
0})\Lambda_1.
\end{align*}
Here $[W(\nu) : V_\lambda]$ denotes the multiplicity 
of the irreducible $U_q({\mathfrak{sl}}(2))$-module $V_\lambda$
with highest weight $\lambda$.
The sum $\sum_{N}^{(\lambda)}$ is the same as (i).
In particular $e^{\Lambda_1}K(\nu)$ is skew-invariant under the Weyl group.
\end{theorem}
\begin{proof}
(i) In view of 
(\ref{eq:Rspecial}),
the two equalities are 
equivalent.
The first one  is due to (\ref{eq:completeness1}) and 
$\mathrm{ch}\,W(\nu) = \prod_{s\ge 1}Q_s^{\nu_s}$ by
{\sc Proposition} \ref{cor:Qcharacter}.
(ii) In view of 
(\ref{eq:Kspecial}),
the two equalities are again equivalent. 
To be self-contained, let us include a 
quick proof of the first one although this has been done in \cite{K}.
Let $\gamma_j$ be as in (\ref{eq:truncatedgamma}) and $\mu \in \Z_{\ge 0}$.
In the expansion of $\prod_{j=1}^l(1-v_j)^{-\beta_j-1}$ by means of
(\ref{eq:basic}),
specialize the variables as
$v_j = 1-\frac{Q_{j-1}Q_{j+1}}{Q^2_j}$ (hence $z_j = Q_1^{-2j}$)
and $\beta_j = \gamma_j - \mu$.
The result reads
\begin{displaymath}
\begin{split}
&Q_1^{-\mu+1}\left(\frac{Q_l}{Q_{l+1}}\right)^{\gamma_l-\mu+1}
\mathrm{ch}\,W(\nu) \\
&\qquad = \sum_N Q_1^{-2\sum_{i=1}^liN_i}
\prod_{j=1}^l
\binom{\gamma_j - \mu + 2\sum_{j < k \le l}(k-j)N_k + N_j}{N_j},
\end{split}
\end{displaymath}
where $\sum_N$ is taken over $N_1,$ $N_2$, \dots $\in
\Z_{\ge 0}$.
Picking up the coefficient of $Q_1^{-\mu}$, we get
\begin{displaymath}
\begin{split}
\sum_{N:\, 2\sum_{i=1}^liN_i = \mu} K_l(\nu,N) &= 
\Res_{Q_1= \infty}\left(Q_1\left(\frac{Q_l}{Q_{l+1}}\right)^{\gamma_l - \mu + 1}
\mathrm{ch}\,W(\nu)\right) \frac{dQ_1}{Q_1}\\
&= \Res_{x = \infty}\left(x(1-x^{-2})
\left(\frac{Q_l}{Q_{l+1}}\right)^{\gamma_l - \mu + 1}
\mathrm{ch}\,W(\nu)\right) \frac{dx}{x},
\end{split}
\end{displaymath}
where $Q_1 = x+x^{-1}$ is used.
In the limit $l \rightarrow \infty$ this is equivalent to 
$K(\nu) = (1-e^{-\alpha_1})\,\mathrm{ch}\,W(\nu)$ due to
(\ref{eq:Kspecial}) and  the property (b) in {\sc
Proposition} \ref{pr:qsys}.
\end{proof}

It is curious that in general 
the sum $\sum_N^{(\lambda)}$ involves the contributions 
{}from those $N$ that do not satisfy the assumption in 
{\sc Theorem} \ref{th:Rnds}.


\section{Discussion}\label{sec:discussion}

In this paper we have proposed the string center
equation (SCE) relevant to the string solutions of the Bethe equation
at $q=0$.
The number of off-diagonal solutions to SCE is identified
with the weight multiplicities of the quantum space by
constructing an explicit combinatorial formula
$R(\nu,N)$.

It is quite common to reduce the Bethe equation to the one
for string centers.
Indeed such analyses have been done extensively at $q=1$,
and has led to the well known fermionic formula
$K(\nu,N)$ \cite{K}.
However at $q=0$, systematic counting of the number of
solutions had been left untouched.
The result in this paper reveals another aspect
of the combinatorial completeness of the
string hypothesis.
The fermionic form $K(\nu,N)$ is relevant to $q=1$
and the multiplicity of irreducible components, while our $R(\nu,N)$ is
relevant to $q=0$ and the weight multiplicities.
Their generating functions are simply related as
(\ref{eq:RKKinf}) and (\ref{eq:RKK2}).

In this paper we have exclusively treated the
$U_q(\selh)$ case.
Many results here admit  straightforward generalizations to
$U_q(X^{(1)}_n)$, which will be the subject of our subsequent paper.
In place of (\ref{eq:Rexpand})--(\ref{eq:Ddef}),
our main formula is
($\nu = (\nu^{(a)}_s), N = (N^{(a)}_m)$)
\begin{gather*}
R(\nu,N) = \sum_{J \subset \N^n} D_J
\prod_{(a,m)\in \N^n \setminus J}
\binom{P^{(a)}_m + N^{(a)}_m}{N^{(a)}_m}
\prod_{(a,m) \in J}
\binom{P^{(a)}_m + N^{(a)}_m - 1}{N^{(a)}_m - 1},
\\
D_J = \begin{cases}
1 & \text{ if } J = \emptyset,\\
\det_{(a,m),(b,k) \in J}
\left((\alpha_a \vert \alpha_b)\min(t_bm, t_ak)-
\delta_{a,b}\delta_{m,k}\right) & \text{otherwise},
\end{cases}
\end{gather*}
where $\N^n = \{(a,m) \mid 1 \le a \le n,
m \in \N \}$, and the other notations are the same as
\cite{HKOTY} under the identification of
$(P^{(a)}_j, N^{(a)}_j)$ here with
$(p^{(a)}_j, m^{(a)}_j)$ there.
As the ${\mathfrak{sl}}(2)$ case, the above $R(\nu,N)$ contains 
the fermionic form in \cite{KR} as 
the summand corresponding to $J = \emptyset$.
With this $R(\nu,N)$,
{\sc Theorem} \ref{th:Rnds}, \ref{th:main},
{\sc Proposition} \ref{pr:qsys} generalize to arbitrary $X_n$.
On the other hand, the identification like
{\sc Proposition} \ref{cor:Qcharacter}, hence
{\sc Theorem} \ref{th:completeness3} are attained
for the  non-exceptional series $X_n = A_n$, $B_n$, $C_n$, $D_n$
only, due to a technical complexity.

Another direction of the generalization is to
seek a $q$-analogue of $R(\nu,N)$ that expresses
the unrestricted one dimensional
configuration sums (1dsums)
over the quantum space $W$ in the sense of
\cite{HKOTY}.
So far we have only obtained a conjecture for the XXZ case
jointly with G.\ Hatayama, M.\ Okado and T.\ Takagi.
%

\vskip0.3cm
{\em Acknowledgment}: The authors thank G.\ Hatayama,
 M.\ Okado and
T.\ Takagi for stimulating discussion and collaboration
on a generalization of the present work.
They also thank M.\ T.\ Batchelor and V.\ O.\ Tarasov for useful 
correspondence.

\appendix

\section{M\"obius function $\mu(\pi, \pi')$}\label{app:mobius}
Let us explain a minimum about the
partition of sets and the M\"obius function on it.
For a more extensive treatment
see \cite{A,B,S}.

\subsection{Partition of set}\label{subapp:1}
Let $N \in \N$.
By definition
$\pi = (\pi_1, \ldots, \pi_l)$ is called a {\em partition\/} of
a set $\{1, \ldots, N\}$ if
\begin{equation*}
\{1, \ldots, N\} = \pi_1 \sqcup \cdots \sqcup \pi_l
\end{equation*}
is a disjoint union decomposition.
Here, the ordering of $\pi_1$, $\pi_2$, \dots, $\pi_l$ does not
matter, e.g.,
$(\pi_1, \pi_2, \ldots, \pi_l)$ and
$(\pi_2, \pi_1, \ldots, \pi_l)$ are the same partition.
Each $\pi_i$ is called a {\em block\/} of $\pi$ and
$l$ is called a {\em length\/} of $\pi$.
Let $L_N$ denote the set of partitions of $\{1, \ldots, N\}$.
Here are the first three:
\begin{align*}
L_1 &= \{ 1 \},\\
L_2 &= \{ 12,\;  1/2\},\\
L_3 &= \{ 123, \; 12/3, \; 13/2, \; 23/1, \; 1/2/3\},
\end{align*}
where, for example, $23/1$ stands for the partition
$\pi = (\pi_1, \pi_2)$ of length $l(\pi) = 2$ consisting of
the blocks $\pi_1=\{2,3\}$ and $\pi_2 = \{1\}$.
\subsection{Poset structure}\label{subapp:2}
One can endow a natural partial order ``$\le $" with the set $L_N$.
Given two partitions $\pi, \pi' \in L_N$, we say
$\pi \le \pi'$ if each block of $\pi'$ is contained in a block of $\pi$.
For $L_2$ in the above we have $12 \le 1/2$, and for $L_3$
\begin{align*}
& 12/3 \\
123 \;\;\le \quad & 13/2 \;\; \le \;\, 1/2/3,\\
&  23/1
\end{align*}
where there is no order among the middle three.
Sometimes $\pi'$ is called a {\em refinement\/} of $\pi$
when $\pi \le \pi'$.
The partition
$\pi_{\text{max}} = 1/2/\cdots/N$
(resp.\ $\pi_{\text{min}} = 12\ldots N$) is the unique
maximal (resp.\ minimal) element in $L_N$ of length
$l(\pi_{\text{max}}) = N$ (resp.\ $l(\pi_{\text{min}}) = 1$).
Clearly the following three axioms hold:
\begin{enumerate}
\item For any $\pi \in L_N$, $\pi \le \pi$. (reflexivity)
\item If $\pi \le \pi'$ and $\pi' \le \pi$, then $\pi = \pi'$. (antisymmetry)
\item If $\pi \le \pi'$ and $\pi' \le \pi''$, then
$\pi \le \pi''$. (transitivity)
\end{enumerate}
Thus $L_N$ equipped with $\le$ is a partially ordered set (poset)
in the sense of \cite{S}.
\subsection{M\"obius function}\label{subapp:3}
Consider an $\vert L_N \vert$ by $\vert L_N \vert$ matrix $\zeta$ defined by
\begin{equation*}
\zeta = \bigl( \zeta(\pi, \pi') \bigr)_{\pi, \pi' \in L_N},\quad
\zeta(\pi, \pi') = \begin{cases}
1 & \text{if } \, \pi \le \pi'\\
0 & \text{otherwise}.
\end{cases}
\end{equation*}
This matrix is upper triangular with all the diagonal elements being $1$.
Thus it has the inverse
\begin{equation*}
\zeta \mu = 1_{L_N}, \quad
\mu = \bigl( \mu(\pi, \pi') \bigr)_{\pi, \pi' \in L_N}.
\end{equation*}
The matrix elements $\mu(\pi, \pi') \in \Z$ are called the M\"obius function of
the poset $L_N$.
Note from  the definition that
$\mu(\pi, \pi) = 1$ for any $\pi$ and $\mu(\pi , \pi') = 0$ unless $\pi \le
\pi'$.
For example in $N = 2$ and $3$ cases in the above they are explicitly given by
\begin{equation}\label{eq:muexample}
\mu = \begin{pmatrix} 1 & -1\\
                        0 & 1
        \end{pmatrix}, \quad
\mu = \begin{pmatrix} 1 & -1 & -1 & -1 & 2 \\
                        0 & 1  &  0 &  0 & -1\\
                        0 & 0  &  1 &  0 & -1\\
                        0 & 0  &  0 &  1 & -1\\
                        0 & 0  &  0 &  0 &  1
        \end{pmatrix}.
\end{equation}
For general $N$, an explicit formula of  $\mu(\pi , \pi')$ is available
\cite{A,B}, but we do not need it in this paper.
\subsection{M\"obius inversion formula}\label{subapp:4}
Given any function $f : L_N \rightarrow \C$,
define another function $g: L_N \rightarrow \C$ by
\begin{equation*}
g(\pi) = \sum_{\pi' \le \pi} f(\pi').
\end{equation*}
This is a composition with  $\zeta$ introduced previously.
In the vector-matrix notation it is expressed as
$g = f \zeta$, hence is equivalent to $f = g \mu$:
\begin{equation*}
f(\pi) = \sum_{\pi' \le \pi} g(\pi')\mu(\pi',\pi),
\end{equation*}
which is the M\"obius inversion formula.
Here is a simple example of its application:
\begin{proposition}\label{pr:inversion1}
Let $X$ be an indeterminate. For any $\pi \in L_N$ we have
\begin{align*}
X^{l(\pi)} &= \sum_{\pi' \le \pi} ( X )_{l(\pi')},\\
( X )_{l(\pi)} &= \sum_{\pi' \le \pi} \mu(\pi',\pi) X^{l(\pi')},
\end{align*}
where $( X )_l = X(X-1) \cdots (X-l+1)$.
\end{proposition}
\begin{proof}
It suffices to show the former assuming that $X$ is
any positive integer.
Notice that
$X^{l(\pi)}$ is the number of maps
$\phi: \{1, \ldots, N\} \rightarrow \{1, \ldots, X\}$ such that
$\phi(i) = \phi(j)$ if $i$ and $j$ belong to the same block of $\pi$.
Similarly, $( X )_{l(\pi)}$ is the number of maps
$\phi: \{1, \ldots, N\} \rightarrow \{1, \ldots, X\}$ such that
$\phi(i) = \phi(j)$ if and only if  $i$ and $j$ belong to the same block of
$\pi$.
Since ``if\/" case consists of the disjoint union of ``if and only
if\/" cases labeled by $\pi'\;( \le \pi)$,
the former relation holds.
\end{proof}
\subsection{Product poset}\label{subapp:5}
Consider the product set
$L_{N_1} \times \cdots \times L_{N_m}$
for any positive integers $N_1$, \dots, $N_m$.
Denote its  elements by
$\pi = (\pi^{(1)}, \ldots, \pi^{(m)})$,
where $\pi^{(i)} \in L_{N_i}$.
One can equip
$L_{N_1} \times \cdots \times L_{N_m}$ with a poset structure by introducing
the partial order as
\begin{displaymath}
(\pi^{(1)},\ldots, \pi^{(m)}) \le
(\pi^{(1)'},\ldots, \pi^{(m)'})
\overset{\mathrm{def}}{\Longleftrightarrow}
\pi^{(i)} \le \pi^{(i)'} \ \text{for any}\ 1 \le i \le m.
\end{displaymath}
The length function is also introduced as
$l\bigl((\pi^{(1)}, \ldots, \pi^{(m)})\bigr) = l(\pi^{(1)}) + \cdots
 + l(\pi^{(m)})$.
The unique maximal element in $L_{N_1} \times \cdots \times L_{N_m}$ is
$\pi_{\text{max}} = (\pi^{(1)}_{\text{max}},\ldots,\pi^{(m)}_{\text{max}})$,
where $\pi^{(i)}_{\text{max}} = 1/2/\cdots/N_i$ is the maximal
one in $L_{N_i}$.
Obviously the M\"obius function
of this poset is the direct product of the one for each component:
\begin{equation*}
\mu\bigl((\pi^{(1)},\ldots,\pi^{(m)}),
(\pi^{(1)'},\ldots, \pi^{(m)'})\bigr) =
\prod_{i=1}^m
\mu_{N_i}(\pi^{(i)},\pi^{(i)'}),
\end{equation*}
where we have written the
M\"obius function of $L_{N_i}$ as $\mu_{N_i}$.
Combining this with {\sc Proposition} \ref{pr:inversion1}, we get
\begin{equation}
\sum_{\pi \in
L_{N_1} \times \cdots \times L_{N_m}} \mu(\pi, \pi_{\text{max}})
X^{l(\pi^{(1)})}_1 \cdots X^{l(\pi^{(m)})}_m
= \prod_{i=1}^m( X_i )_{l(\pi^{(i)}_{\text{max}})}
= \prod_{i=1}^m( X_i )_{N_i},\label{eq:musum}
\end{equation}
where $X_1$, \dots, $X_m$ are indeterminates.
In the main text we use the M\"obius inversion formula for the poset
$L_{N_1} \times \cdots \times L_{N_m}$ and (\ref{eq:musum}).


\end{document}